\documentclass[12pt]{article}
\usepackage{amssymb,a4}
\usepackage{mathrsfs}
\usepackage{epsf}
\usepackage{bm}
\usepackage{latexsym}
\usepackage{amsmath,amsfonts,amssymb,amsthm}
\RequirePackage{amsopn} \RequirePackage{amsfonts}
\RequirePackage{amsthm}
\newcommand{\f}{\mathfrak}

\newcommand{\al}{\alpha}

\newcommand{\charge}[1]{V_{L}^{#1}}

\newcommand{\charlam}[1]{V_{#1+L}}

\newcommand{\1}{{{\bf 1}}}
\newcommand{\id}{{\rm id}}

\newcommand{\End}{{\rm End}\,}

\newcommand{\Z}{\mathbb{Z}}
\newcommand{\Q}{\mathbb{Q}}
\newcommand{\C}{\mathbb{C}}

\def\wt{{\rm wt}}
\def\De{\Delta}
\def\de{\delta}
\def\be{\beta}

\newcommand{\la}{\lambda}\newcommand{\La}{\Lambda}

\def\C{{\mathbb C}}

\def\Z{{\mathbb Z}}

\def\1{{\bf 1}}
\def \End{{\rm End}}

\def \pf{\noindent {\bf Proof: \,}}
\def\theequation{5.\arabic{equation}}
\def \a{\alpha}
\def \b{\beta}
\def \e{\epsilon}
\def \h{\mathfrak{h}}
\def \l{\lambda}
\def \w{\omega}
\begin{document}
\newtheorem{theorem}{Theorem}[section]
\newtheorem{proposition}[theorem]{Proposition}
\newtheorem{lem}[theorem]{Lemma}
\newtheorem{corollary}[theorem]{Corollary}
\newtheorem{definition}[theorem]{Definition}
\newtheorem{example}[theorem]{Example}
\newtheorem{remark}[theorem]{Remark}

\numberwithin{equation}{section}

\newenvironment{namelist}[1]{%
\begin{list}{}
{
\let\makelabel\namelistlabel
\settowidth{\labelwidth}{#1}
 \setlength{\leftmargin}{1.1\labelwidth}
 }
 }{%
 \end{list}}
 \begin{center}
{\Large {\bf Rationality of vertex operator algebra $V_L^+:$ higher rank}} \\

\vspace{0.5cm} Chongying Dong\footnote{Supported by NSF grants and  a Faculty research grant from  the
University of California at Santa Cruz.}
\\
School of Mathematics,  Sichuan University, Chengdu, 610065 China\\
\& Department of Mathematics, University of
California\\ Santa Cruz, CA 95064 \\
\vspace{.2 cm}
Cuipo Jiang\footnote{Supported  by China NSF grants
10931006,10871125, the Innovation Program of Shanghai Municipal
Education Commission (11ZZ18) and a grant of Science and Technology
Commission
of Shanghai Municipality (No. 09XD1402500).}\\
 Department of Mathematics, Shanghai Jiaotong University\\
Shanghai 200240 China\\
\vspace{.2 cm}
Xingjun Lin\\
Department of Mathematics, Sichuan University\\
 Chengdu, 610065 China
\end{center}

\begin{abstract}
The lattice vertex operator $V_L$ associated to a positive definite
even lattice $L$ has an automorphism of order 2 lifted from $-1$
isometry of $L$. It is established that the fixed point vertex
operator algebra $V_L^+$ is rational.
\end{abstract}
\section{Introduction \label{intro}}
\def\theequation{1.\arabic{equation}}
\setcounter{equation}{0}

The notion of rational vertex operator algebra is an analogue to
that of semisimple Lie algebra or semisimple associative algebra.
Rational vertex operator algebras whose admissible module category
is semisimple form a fundamental class of vertex operator algebras.
Familiar examples of rational vertex operator algebras include the
vertex operator algebras $V_L$ associated with even lattices [D1],
\cite{DLM1}, vertex operator algebras  associated to the irreducible
vacuum representations for affine Kac-Moody algebras with positive
levels [FZ], \cite{DL}, \cite{LL}, vertex operator algebras
associated to the minimal series for the Virasoro algebra \cite{W}.

Let $V$ be a vertex operator algebra and $G$  a finite automorphism
group of $V,$ then the space of $G$-invariants $V^G$ is itself a
vertex operator algebra. A well known conjecture in the orbifold
conformal field theory states that if $V$ is rational  then $V^G$ is
rational. Solving this conjecture has significant applications in
the theory of vertex operator algebras and conformal field theory.

Let $V_L$ be a lattice vertex operator algebra associated with a
positive definite even lattice $L$ which plays an important role in
shaping the theory of vertex operator algebras. The vertex operator
algebra $V_L$ has an automorphism $\theta$ of order $2$ lifted from
the $-1$ isometry of $L$ and we denote the $\theta$-fixed points of
$V_L$ by $V_L^+.$ The vertex operator algebras $V_L^+$ have been
studied extensively. The irreducible modules for $V_L^+$ were
classified in \cite{DN2} and \cite{AD}. The fusion rules among
irreducible $V_{L}^+$-modules were computed in \cite{ADL}. The
$C_2$-cofiniteness of $V_L^+$ was obtained in \cite{Y} and
\cite{ABD}. But the rationality of $V_L^+$ has only been established
if $L$ has rank one \cite{A2}, or $L$ is a special lattice
\cite{DGH}. In this paper, we extend the rationality result to any
lattice. That is, $V_L^+$ is rational for any positive definite even
lattice $L.$

We prove the rationality of $V_L^+$ in three steps. First, we show
that the Zhu algebra $A(V_L^+)$ is a finite dimensional semisimple
associative algebra. Although $A(V_L^+)$ was studied in \cite{DN2}
and \cite{AD} in great detail for the purpose of classification of
irreducible modules, there is still a distance to claim the
semi-simplicity. Second, we use the rationality of $V_L^+$ in the
case that $L$ has rank one and the fusion rules to deal with the
rationality of $V_L^+$ if $L$ has an orthogonal base. Last, we use
the rationality of $V_{L_1}^+$ to prove the rationality of $V_L^+$
for any $L$ where $L_1$ is a sublattice of $L$ with the same rank
and has an orthogonal base. The main idea in the proof of
rationality is to prove that if there is a $V_L^+$-module exact
sequence
$$0\to M^1\to M\to M^2\to 0$$
for any irreducible $V_L^+$-modules $M^1$ and $M^2,$ then  $M$ is a direct sum of $M^1$ and $M^2.$

The representation theory of $V_L^+$ is complete in some sense after
this paper. But the structure theory of $V_L^+$ is far from over.
Determining the derivation Lie algebra and the automorphism group of
$V_L^+$ for an arbitrary positive definite even $L$ remains a major
problem. This has been achieved when the rank of $L$ is one, two
or three \cite{DG1}, \cite{DG2}, \cite{S1}, \cite{S2}, or $L$ is 
a special lattice which is either unimodular or does not have 
roots \cite{S1}, \cite{S2}. Extending these results to general lattice seems
a big challenge.

The paper is organized as follows: In Section 2, we recall
definitions of vertex operator algebra, module, intertwining
operator and fusion rules. We also give some basic facts on vertex
operator algebras in this section. In Section 3, we present the
results on vertex operator algebras $M(1)^+$ and $V_L^+$ and the
irreducible modules. Section 4 is devoted to the proof of the
semi-simplicity of $A(V_L^+).$ We deal with the rationality of
$V_L^+$ when $L$ has an orthogonal base in Section 5. The
rationality of $V_L^+$ for any positive definite even lattice $L$ is
given in Section 6.

\section{Preliminaries}
\def\theequation{2.\arabic{equation}}
\setcounter{equation}{0}

In this section we briefly review the definitions of twisted modules and rationality
from \cite{FLM}, \cite{DLM2}. We present the Zhu algebra, the tensor product
 and fusion rules from \cite{FHL} and \cite{Z}. We also discuss the extensions of modules following \cite{A2} and
 give a sufficient condition under which the extensions are trivial. This result will be used extensively in latter sections.

A {\em vertex operator algebra} $V$ is a $\Z$-graded vector space
$V=\bigoplus_{n\in \Z} V_n$ equipped with a linear map $Y:V\to(\End
V)[[z,z^{-1}]],\,a\mapsto Y(a,z)=\sum_{n\in \Z}a_nz^{-n-1}$ for
$a\in V$ such that $\dim V_n$ is finite for all integer $n$ and that
$V_n=0$ for sufficiently small  $n$ (see \cite{FLM}). There are two
distinguished vectors, the {\it vacuum vector} $\1\in V_0$ and the
{\it Virasoro element} $\w\in V_2$.  By definition $Y(\1,z)={\rm
id}_{V}$, and the component operators $\{L(n)|n\in\Z\}$ of
$Y(\w,z)=\sum_{n\in\Z}L(n)z^{-n-2}$ give a representation of the
Virasoro algebra on $V$ with central charge $c.$ Each homogeneous
subspace $V_n\,(n\geq0)$ is an eigenspace for $L(0)$ with eigenvalue
$n$.

An automorphism $g$ of a vertex operator algebra $V$ is a linear
isomorphism of $V$ satisfying $g(\w)=\w$ and
$gY(a,z)g^{-1}=Y(g(a),z)$ for any $a\in V$. We denote by ${\rm
Aut}(V)$ the group of all automorphisms of $V$. For a subgroup
$G<{\rm Aut}(V)$, the fixed point set $V^G=\{a\in
V\,|\,g(a)=a,\,g\in G\,\}$ has a canonical vertex operator algebra
structure.

Let $g$ be an automorphism of a vertex operator algebra $V$ of order
$T$.  Then $V$ is a direct sum of eigenspaces for $g$:
\[
V=\bigoplus_{r=0}^{T-1}V^{r},\,V^{r}=\{\,a\in V\,|\,g(a)=e^{-\frac{2
\pi ir}{T}} a\,\}.
\]
\begin{definition}{\rm
A {\em weak $g$-twisted $V$-module} $M$ is a vector space equipped
with a linear map
\begin{align*}
Y_{M}:V&\to (\End M)\{z\},\\
a&\mapsto Y_{M}(a,z)=\sum_{n\in\Q}a_nz^{-n-1},\,a_n\in \End M
\end{align*}
such that the following conditions hold for $0\leq r\leq T-1,\,a\in
V^{r},\,b\in V$ and $u\in M$:

(1) $b_mu=0$ if $m$ is sufficiently large,

(2) $Y_{M}(a,z)=\sum_{n\in\Z+\frac{r}{T}}a_nz^{-n-1}$,

(3) $Y_{M}(\1,z)=\id _{M}$,

(4) (the twisted Jacobi identity)
\begin{align*}
\begin{split}
&z_{0}^{-1}\delta\left(\frac{z_{1}-z_{2}}{z_{0}}\right)Y_{M}(a,z_{1})Y_{M}(b,z_{2})-z_{0}^{-1}\delta\left(
\frac{z_{2}-z_{1}}{-z_{0}}\right)Y_{M}(b,z_{2})Y_{M}(a,z_{1})\\
&\quad=z_{2}^{-1}\left(\frac{z_{1}-z_{0}}{z_{2}}\right)^{-\frac{r}{T}}\delta\left(\frac{z_{1}-z_{0}}{z_{2}}\right)Y_{M}(Y(a,z_{0})b,z_{2}).
\end{split}
\end{align*}
}
\end{definition}

A weak $g$-twisted $V$-module is denoted by $(M,\,Y_{M})$, or simply
by $M$. In the case $g$ is the identity, any weak $g$-twisted
$V$-module is called a \textit{weak $V$-module}. A
\textit{$g$-twisted weak $V$-submodule} of a $g$-twisted weak module
$M$ is a subspace $N$ of $M$ such that $a_nN\subset N$ hold for all
$a\in V$ and $n\in\Q$. If $M$ has no $g$-twisted weak $V$-submodule
except $0$ and $M$, $M$ is called \textit{irreducible} or
\textit{simple}.

Set $Y_{M}(\w,z)=\sum_{n\in\Z}L(n)z^{-n-2}$. Then
$\{L(n)\,|\,n\in\Z\,\}$ give a representation of the Virasoro
algebra on $M$ with central charge $c$ and the $L(-1)$-derivative
property
\begin{equation}\label{DP1}
Y_{M}(L(-1)a,z)=\frac{d}{dz}Y_M(a,z)\hbox{ for all }a\in V
\end{equation}
holds for any $a\in V$ (see \cite{DLM1}).

\begin{definition}{\rm
An \textit{admissible $g$-twisted $V$-module} $M$ is a weak
$g$-twisted $V$-module which has a $\frac{1}{T}\Z_{\geq
0}$-gradation $M=\bigoplus_{n\in\frac{1}{T}\Z_{\geq 0}}M(n)$ such
that
\begin{align}\label{AD1}
a_mM(n)\subset M(\wt{a}+n-m-1)
\end{align}
for any homogeneous $a\in V$ and $m,\,n\in\Q$. }
\end{definition}
In the case $g$ is the identity, any admissible $g$-twisted
$V$-module is called an {\em admissible $V$-module}. Any $g$-twisted
weak $V$-submodule $N$ of a $g$-twisted admissible $V$-module is
called a {\it $g$-twisted admissible $V$-submodule} if
$N=\bigoplus_{n\in \frac{1}{T}\Z_{\geq 0}} N\cap M(n)$.

A $g$-twisted admissible $V$-module $M$ is said to be
\textit{irreducible} if $M$ has no non-trivial admissible weak
$V$-submodule. When a $g$-twisted admissible $V$-module $M$ is a
direct sum of irreducible admissible submodules, $M$ is called
\textit{completely reducible}.

\begin{definition}
{\rm
A vertex operator algebra $V$ is said to be \textit{$g$-rational} if
any $g$-twisted admissible $V$-module is completely reducible. If
$V$ is $\id_{V}$-rational, then $V$ is called \textit{rational}.}
\end{definition}

\begin{definition}{\rm
A \textit{$g$-twisted $V$-module} $M=\bigoplus_{\lambda\in\C}
M_{\lambda}$ is a $\C$-graded weak $g$-twisted $V$-module with
$M_\lambda=\{u\in M|L(0)u=\lambda u\}$ such that $M_\lambda$ is
finite dimensional and for fixed $\lambda\in\C$, $M_{\lambda+n/T}=0$
for sufficiently small integer $n$. A vector $w\in M_\l$ is called a
weight vector of weight $\l,$ and we write $\l=\wt{w}.$ }
\end{definition}

We next define  Zhu's algebra $A(V)$ which is an associative algebra
following [Z]. For any homogeneous vectors $a, b\in V$, we define
$$a\ast b= Res_z\frac{(1+z)^{wt a}}{z}Y(a, z)b,$$ $$a\circ
b= Res_z\frac{(1+z)^{wt a}}{z^2}Y(a, z)b$$ and extend the two
operations to $V\times V$ bilinearly. Denote by $O(V)$ the subspace
of $V$ linearly spanned  by  $a\circ b$, for all $a,b\in V$ and set
$A(V)=V/O(V)$.   The following theorem is due to \cite{Z}.
\begin{theorem}\label{P3.1}
(1) The bilinear operation $*$ induces $A(V)$ an associative algebra
structure. The vector $[\1]$ is the identity and $[\w]$ is in the
center of $A(V)$.

(2) Let $M=\bigoplus_{n=0}^{\infty}M(n)$ be an admissible $V$-module
with $M(0)\ne 0.$ Then the linear map
\[
o:V\rightarrow\End M(0),\;a\mapsto o(a)|_{M(0)}
\]
induces an algebra homomorphism from $A(V)$ to $\End M(0)$. Thus
$M(0)$ is a left $A(V)$-module.

(3) The map $M\mapsto M(0)$ induces a bijection between the set of
equivalence classes of irreducible admissible $V$-modules and the
set of equivalence classes of irreducible $A(V)$-modules.
\end{theorem}

Now we consider the tensor product vertex algebra and the tensor
product modules for tensor product vertex operator algebra. The
tensor product of vertex operator algebras $(V^1, Y, \1, \omega^1), \cdots, (V^p, Y, \1,\omega^p)$ is
constructed on the tensor product vector space$$V = V^1\otimes
\cdots \otimes V^p$$
where the vertex operator $Y(\cdot,z)$ is
defined by
$$Y(v^1\otimes \cdots \otimes v^p, z) = Y(v^1, z)\otimes
\cdots \otimes Y(v^p, z)$$\\for $v^i \in V^i\ (1\leq i \leq p),$
the vacuum vector is $$1 = 1\otimes \cdots \otimes 1$$
and the Virasoro
element is
$$\omega = \omega^1 \otimes \cdots \otimes 1 + \cdots +1\otimes
\cdots\otimes \omega^p.$$
then  $(V, Y, 1, \omega)$
is a vertex operator algebra (see  \cite{FHL}, \cite{LL}).

Let $M^i$  be an admissible $V^i$-module for $i=1,...,p.$  We may construct
the tensor product admissible module $M^1\otimes \cdots \otimes M^p$ for the
tensor product vertex operator algebra $V^1\otimes \cdots \otimes
V^p$ by
$$Y(v^1\otimes \cdots \otimes v^p) =
Y(v^1, z)\otimes \cdots \otimes Y(v^p, z).$$ Then  $(M^1\otimes \cdots \otimes M^p, Y)$ is an admissible $V^1\otimes
\cdots \otimes V^p$-module. The following result was given
in \cite{FHL} and \cite{DMZ}.
\begin{theorem}\label{t2.6}
Let $V^1, \cdots, V^p$ be rational vertex operator algebras, then
$V^1\otimes \cdots \otimes V^p$  is rational and any irreducible $V^1\otimes \cdots \otimes V^p$-module is a tensor
product $M^1\otimes \cdots \otimes M^p$ , where $M^1, \cdots, M^p$ are
some irreducible modules for the vertex operator algebras $V^1,
\cdots, V^p,$ respectively.
\end{theorem}

Let $M = \bigoplus_{r\in \mathbb{C}}{M(r)}$ be a $V$-module. Set $M'
= \bigoplus_{\lambda \in \mathbb{C}}{M_\lambda^*}$, the restricted
dual of $M$. It was proved in [FHL] that $M'$ is naturally a
$V$-module where the vertex operator map denoted by $Y'$ is defined
by the property
$$\langle Y'(a, z)u', v\rangle  = \langle u', Y(e^{zL(1)}(-z^{-2})^{L(0)}a, z^{-1})v\rangle $$\\for $a\in V, u'\in
M'$ and $v\in M$. The $V$-module $M'$ is called the contragredient
module of $M$. It was proved that if $M$ is irreducible ,  then so
is $M'$. A $V$-module $M$ is said to be self dual if $M$ and $M'$
are isomorphic $V$-modules.

We now recall the  notion of intertwining operators and fusion
rules from [FHL].
\begin{definition}
Let $M^1$, $M^2$, $M^3$ be weak $V$-modules. An intertwining
operator $\mathcal {Y}( \cdot , z)$ of type $\left(\begin{tabular}{c}
$M^3$\\
$M^1$ $M^2$\\
\end{tabular}\right)$ is a linear map$$\mathcal
{Y}(\cdot, z): M^1\rightarrow Hom(M^2, M^3)\{z\}$$ $$v^1\mapsto
\mathcal {Y}(v^1, z) = \sum_{n\in \mathbb{C}}{v_n^1z^{-n-1}}$$
satisfying the following conditions:

(1) For any $v^1\in M^1, v^2\in M^2$and $\lambda \in \mathbb{C},
v_{n+\lambda}^1v^2 = 0$ for $n\in \mathbb{Z}$ sufficiently large.

(2) For any $a \in V, v^1\in M^1$,
$$z_0^{-1}\delta(\frac{z_1-z_2}{z_0})Y_{M^3}(a, z_1)\mathcal
{Y}(v^1, z_2)-z_0^{-1}\delta(\frac{z_1-z_2}{-z_0})\mathcal{Y}(v^1,
z_2)Y_{M^2}(a, z_1)$$
$$=z_2^{-1}\delta(\frac{z_1-z_0}{z_2})\mathcal{Y}(Y_{M^1}(a, z_0)v^1, z_2).$$

(3) For $v^1\in M^1$, $\dfrac{d}{dz}\mathcal{Y}(v^1,
z)=\mathcal{Y}(L(-1)v^1, z)$.
\end{definition}
All of the intertwining operators of type $\left(\begin{tabular}{c}
$M^3$\\
$M^1$ $M^2$\\
\end{tabular}\right)$ form a vector space denoted by $I_V\left(\begin{tabular}{c}
$M^3$\\
$M^1$ $M^2$\\
\end{tabular}\right)$. The dimension of $I_V\left(\begin{tabular}{c}
$M^3$\\
$M^1$ $M^2$\\
\end{tabular}\right)$ is called the
fusion rule of type $\left(\begin{tabular}{c}
$M^3$\\
$M^1$ $M^2$\\
\end{tabular}\right)$ for $V$.

We now have the following result which was essentially proved in
[ADL].

\begin{theorem}\label{n2.2c}
Let $V^1, V^2$ be rational vertex operator algebras. Let $M^1 , M^2,
M^3$ be $V^1$-modules and  $N^1, N^2, N^3$ be $V^2$-modules such
that
$$dim I_{V^1}\left(\begin{tabular}{c}
$M^3$\\
$M^1$ $M^2$\\
\end{tabular}\right)< \infty , \ dim I_{V^2}\left(\begin{tabular}{c}
$N^3$\\
$N^1$ $N^2$\\
\end{tabular}\right)< \infty.$$
Then the linear map
$$\sigma: I_{V^1}\left(\begin{tabular}{c}
$M^3$\\
$M^1$ $M^2$\\
\end{tabular}\right)\otimes I_{V^2}\left(\begin{tabular}{c}
$N^3$\\
$N^1$ $N^2$\\
\end{tabular}\right)\rightarrow I_{V^1\otimes V^2}\left( \begin{tabular}{c}
$M^3\otimes N^3$\\
$M^1\otimes N^1$ $M^2\otimes N^2$\\
\end{tabular}\right)$$
$$\mathcal{Y}_1( \cdot, z)\otimes \mathcal{Y}_2(\cdot, z)\mapsto (\mathcal{Y}_1\otimes
\mathcal{Y}_2)(\cdot, z)$$\\is an isomorphism, where
$(\mathcal{Y}_1\otimes \mathcal{Y}_2)(\cdot, z)$ is defined by
$$(\mathcal{Y}_1\otimes \mathcal{Y}_2)(\cdot, z)(u^1\otimes v^1,z)u^2\otimes v^2 = \mathcal{Y}_1(u^1,z)u^2\otimes \mathcal{Y}_2(v^1,z)v^2.$$
\end{theorem}

Now let $M^1, M^2$ be weak $V$-modules, we call a weak $V$-module
$M$ an extension of $M^2$ by $M^1$ if there is a short exact
sequence
$$0\rightarrow M^1\rightarrow M \rightarrow M^2\rightarrow 0.$$
Then we could define the equivalence of two extensions, and then
define the extension group $Ext_V^1(M^2, M^1)$. Recall that $V$ is
called $C_2$-cofinite if the subspace $C_2(V)$ of $V$ spanned by
$u_{-2}v$ for $u,v\in V$ has finite codimension. We have the
following facts \cite{A1}.

\begin{theorem}\label{n2.1}
Let M and N be irreducible V-modules, then $Ext_V^1(N, M) = 0$ if
and only if $Ext_V^1(M' , N') = 0.$
\end{theorem}

\begin{theorem}\label{n2.2}
Let $V$ be a $C_2$-cofinite vertex operator algebra, then $V$ is
rational if and only if $Ext_V^1(N, M) = 0$ for any pair of
irreducible $V$-modules $M$ and $N$. \end{theorem}

The following result will be extensively used in Sections 5 and 6.
\begin{lem}\label{la1}
Let $V$ be a vertex operator algebra and $U$ a rational vertex operator subalgebra
of $V$ with the same Virasoro element. Let $M^1,M^2$ be irreducible $V$-modules.
Assume that
$$I_{U}\left(\begin{tabular}{c}
$N^1$\\
$N$\  $N^2$
\end{tabular}\right)=0$$
for any irreducible $U$-submodules $N^1,N,N^2$ of $M^1,V,M^2,$ respectively.
Then $Ext_V^1(M^2, M^1) = 0.$
\end{lem}

\pf Let $M$ be an extension of $M^2$ by $M^1.$ Then
$M=M^1\oplus M^2$ as $U$-modules as $U$ is rational. Let $N^1,N,N^2$
be any irreducible $U$-submodules of $M^1,V,M^2,$ respectively. Then
$P_{N^1}Y(u,z)|_{N^2}$ for $u\in N$ is an intertwining operator of
type $I_{U}\left(\begin{tabular}{c}
$N^1$\\
$N$\  $N^2$
\end{tabular}\right)$ where $P_{N^1}$ is the projection from $M$ to $N_1.$ From the assumption, $P_{N^1}Y(u,z)|_{N^2}=0.$ Since $N^1, N, N^2$ are arbitrary, we see that
$u_nM^2\subset M^2$ for any $u\in V$ and $n\in \Z.$ As a result, $M^2$ is a $V$-module
and $Ext_V^1(M^2, M^1) = 0.$ The proof is complete.
\qed

\section{Vertex operator algebras $M(1)^+$ and $V_L^+$ }
In this section we recall vertex operator algebras $M(1)^+$ and
$V_L^+$ \cite{FLM} and related results \cite{DN1}, \cite{DN2},
\cite{DN3}, \cite{A1}, \cite{A2}, \cite{AD}, \cite{ADL}.  In particular,
the irreducible modules, fusions and contragredient modules of irreducible modules
for $V_L^+$ are discussed.

Let $L$ be a positive definite even lattice in the sense that $L$ has a
symmetric positive definite  $\mathbb{Z}$-valued
$\mathbb{Z}$-bilinear form $(\cdot, \cdot)$ such that $(\alpha,
\alpha)\in 2\mathbb{Z}$ for any $\alpha \in L$.
We set $\h=\C\otimes_{\Z} L$ and extend $(\cdot\,,\cdot)$ to a
$\C$-bilinear form on $\h$. Let
$\hat{\h}=\C[t,t^{-1}]\otimes\h\oplus\C C$ be the affinization of
commutative Lie algebra $\h$ defined by
\begin{align*}
[\beta_1\otimes t^{m},\,\beta_2\otimes
t^{n}]=m(\beta_1,\beta_2)\delta_{m,-n}C\hbox{ and }[C,\hat{\h}]=0
\end{align*}
for any $\beta_i\in\h,\,m,\,n\in\Z$. Then $\hat{\h}^{\geq
0}=\C[t]\otimes\h\oplus\C C$ is a commutative subalgebra. For any
$\lambda\in\h$, we can define a one dimensional $\hat{\h}^{\geq
0}$-module $\C e^\lambda$ by the actions $\rho(h\otimes
t^{m})e^\lambda=(\lambda,h)\delta_{m,0}e^\lambda$ and
$\rho(C)e^\lambda=e^\lambda$ for $h\in\h$ and $m\geq0$. Now we
denote by
\begin{align*}
M(1,{\lambda})=U(\hat{\h})\otimes_{U(\hat{\h}^{\geq 0})}\C
e^\lambda\cong S(t^{-1}\C[t^{-1}]),
\end{align*}
the $\hat{\h}$-module induced from $\hat{\h}^{\geq 0}$-module $\C
e^\lambda$. Set $M(1)=M(1,0).$ Then there exists a linear map
$Y:M(1)\to(\End M(1,\lambda)[[z,z^{-1}]]$ such that
$(M(1),\,Y,\,\1,\,\w)$ has a simple vertex operator algebra
structure and $(M(1,\lambda),Y)$ becomes an irreducible
$M(1)$-module for any $\lambda\in\h$ (see \cite{FLM}). The vacuum
vector and the Virasoro element are given by $\1=e^0$ and
$\w=\frac{1}{2}\sum_{a=1}^{d}h_a(-1)^2\1$ respectively, where
$\{h_a\}$ is an orthonormal basis of $\h$.

Let $\widehat{L}$ be the canonical central extension of $L$ by $\langle
\kappa\rangle= \langle\kappa|\kappa^2= 1\rangle:$
$$1\rightarrow \langle\kappa
\rangle
\rightarrow\widehat{L}\stackrel {-}{\rightarrow} L\rightarrow 1$$
with the commutator map $c(\al,\be)=\kappa^{(\al,\be)}$ for
$\al,\be\in L$. Let $e: L\rightarrow\widehat{L}$ be a section such
that $e_{0}=1$ and $\epsilon: L\times L\rightarrow\langle \kappa\rangle $  the
corresponding 2-cocycle. We may assume that $\epsilon$ is
bimultiplicative. Then
$\epsilon(\al,\be)\epsilon(\be,\al)=\kappa^{(\al,\be)}$,
$$
\epsilon(\al,\be)\epsilon(\al+\be,\gamma)=\epsilon(\be,\gamma)(\al,\be+\gamma),$$
and $e_{\al}e_{\be}=\epsilon(\al,\be)e_{\al+\be}$ for
$\al,\be,\gamma\in L$. Let $\theta$ denote the automorphism of
$\widehat{L}$ defined by $\theta(e_{\al})=e_{-\al}$ and
$\theta(\kappa)=\kappa$. Set
$K=\{a^{-1}\theta(a)|a\in\widehat{L}\}$. Note that if $(\alpha,\beta)\in 2\Z$ for
all $\alpha, \beta\in L$ then the $\widehat{L}=L\times \langle\kappa
\rangle$ is a direct product of abelian groups. In this case we can and do choose
 $\epsilon(\al,\be)=1$ for all $\alpha,\beta\in L.$

The lattice vertex
operator algebra associated to $L$ is given by
$$V_L=M(1)\otimes \C^{\epsilon}[L],$$
 where
$\C^{\epsilon}[L]$ is the twisted group algebra of $L$ with a basis $e^{\alpha}$ for
$\alpha\in L$ and is an $\widehat L$-module such that $e_\alpha
e^\beta= \epsilon(\al,\be)e^{\alpha+\beta}.$ Note that if $(\al,\be)\in 2\Z$ for
all $\alpha,\beta\in L$ then $\C^{\epsilon}[L]=\C[L]$ is the usual group algebra.

Recall that
$L^{\circ}=\{\,\lambda\in\h\,|\,(\alpha,\lambda)\in\Z\,\}$
is the dual lattice of $L$.  There is an $\widehat{L}$-module
structure on $\C[L^{\circ}]=\bigoplus_{\lambda\in
L^{\circ}}\C e^\lambda$ such that $\kappa$ acts as $-1$
(see \cite{DL}). Let $L^{\circ}=\cup_{i\in
L^{\circ}/L}(L+\lambda_i)$ be the coset decomposition such that
$(\lambda_i,\lambda_i)$ is minimal among all $(\lambda,\lambda)$ for
$\lambda\in L+\lambda_i.$ In particular, $\lambda_0=0.$
Set $\C[L+\lambda_i]=\bigoplus_{\alpha\in L}\C
e^{\alpha+\lambda_i}.$ Then $\C[L^{\circ}]=\bigoplus_{i\in
L^{\circ}/L}\C[L+\lambda_i]$ and each $\C[L+\lambda_i]$ is an
$\widehat L$-submodule of $\C[L^{\circ}].$ The action of
$\widehat L$ on $\C[L+\lambda_i]$ is as follows:
$$e_{\alpha}e^{\beta+\lambda_i}=\e(\a,\b)e^{\a+\b+\l_i}$$
for $\alpha,\,\b\in L.$ On the surface, the module structure on each
$\C[L+\lambda_i]$ depends on the choice of $\lambda_i$ in
$L+\lambda_i.$ It is easy to prove that different choices of
$\lambda_i$ give isomorphic $\widehat L$-modules.

 Set $\C[M]=\bigoplus_{\lambda\in M}\C e^{\lambda}$ for a subset $M$ of
$L^{\circ}$, and define $V_M=M(1)\otimes\C[M]$. Then $V_L$ is a
rational vertex operator algebra and $V_{L+\lambda_i}$ for $i\in
L^{\circ}/L$ are the irreducible modules for $V_L$ (see \cite{B},
\cite{FLM}, \cite{D1}, \cite{DLM1}).

Define a linear isomorphism $\theta:V_{L+\lambda_i}\to
V_{L-\lambda_i}$ for $i\in L^{\circ}/L$ by
\begin{align*}
\theta(\beta_{1}(-n_{1})\beta_{2}(-n_{2})\cdots \beta_{k}(-n_{k})
e^{\alpha+\lambda_i})=(-1)^{k}\beta_{1}(-n_{1})\beta_{2}(-n_{2})\cdots
\beta_{k}(-n_{k}) e^{-\alpha-\lambda_i}
\end{align*}
for $\beta_i\in\h,\,n_i\geq1$ and $\alpha\in L$ if
$2\lambda_i\not\in L,$ and
\begin{multline*}
\theta(\beta_{1}(-n_{1})\beta_{2}(-n_{2})\cdots \beta_{k}(-n_{k})e^{\alpha+\lambda_i})\\
=(-1)^{k}c_{2\l_i}\e(\alpha,2\lambda_i)\beta_{1}(-n_{1})\beta_{2}(-n_{2})\cdots
\beta_{k}(-n_{k})e^{-\alpha-\lambda_i}
\end{multline*}
if $2\lambda_i\in L$ where $c_{2\l_i}$ is a square root of
$\e(2\l_i,2\l_i).$
  Then $\theta$ defines a linear
isomorphism from $V_{L^{\circ}}$ to itself such that
 $$\theta Y(u,z)v=Y(\theta u,z)\theta v$$
for $u\in V_{L}$ and $v\in V_{L^{\circ}}.$ In particular, $\theta$ is
an automorphism of $V_{L}$ which induces an automorphism of $M(1).$

For any $\theta$-stable subspace $U$ of $V_{L^{\circ}}$, let $U^\pm$
be the $\pm1$-eigenspace of $U$ for $\theta$. Then $V_L^+$ is a
simple vertex operator algebra.

Also recall the $\theta$-twisted Heisenberg algebra $\h[-1]$ and its
irreducible module $M(1)(\theta)$ from \cite{FLM}. Let
$\chi$ be a central character of $\widehat{L}/K$ such that
$\chi(\kappa)=-1$ and $T_{\chi}$ the irreducible
$\widehat{L}/K$-module with central character $\chi$. Note that if
$(\al,\be)\in2\Z$ for all $\al,\be\in L$ then
$\widehat{L}/K=L/K\times \langle \kappa\rangle$ and each $T_{\chi}$
is, in fact, an $L/2L$-module. In particular, $T_{\chi}$ is
one-dimensional. In this case, let $\be_1,\cdots,\be_d$ be a basis
of $L$, then $T_{\chi}=T_{\chi_1}\otimes \cdots \otimes T_{\chi_d}$
where each $T_{\chi_i}$ is an irreducible
$\Z\beta_i/2\Z\beta_i$-module such that $e_{\beta_i}$ acts as
$\chi(e_{\beta_i}).$ This fact will be used later.

It is well known that $V_L^{T_{\chi}}=M(1)(\theta)\otimes T_{\chi}$ is an irreducible $\theta$-twisted $V_L$-module
(see \cite{FLM}, \cite{D2}). We  define  actions of
$\theta$ on  $M(1)(\theta)$ and $V_L^{T_{\chi}}$ by
\begin{align*}
\theta(\beta_{1}(-n_{1})\beta_{2}(-n_{2})\cdots
\beta_{k}(-n_{k}))=(-1)^{k}\beta_{1}(-n_{1})\beta_{2}(-n_{2})\cdots
\beta_{k}(-n_{k})
\end{align*}
\begin{align*}
\theta(\beta_{1}(-n_{1})\beta_{2}(-n_{2})\cdots
\beta_{k}(-n_{k})t)=(-1)^{k}\beta_{1}(-n_{1})\beta_{2}(-n_{2})\cdots
\beta_{k}(-n_{k})t
\end{align*}
for $\beta_i\in\h,\,n_i\in \frac{1}{2}+\Z_{+}$ and $t\in T_{\chi}$.
We denote the $\pm 1$-eigenspace of $M(1)(\theta)$ and
$V_L^{T_{\chi}}$ under $\theta$ by $M(1)(\theta)^{\pm}$ and
$(V_L^{T_{\chi}})^{\pm}$ respectively. We have the following
results proved in \cite{DN1}, \cite{DN3} and \cite{AD}:
\begin{theorem}\label{t32}
Any irreducible module for the vertex operator algebra $M(1)^+$ is
isomorphic to one of the following modules:$$ M(1)^+, M(1)^-, M(1,
\lambda) \cong M(1, -\lambda)\ (0\neq \lambda \in \h), M(1)(\theta)^+,
M(1)(\theta)^- .$$
\end{theorem}
\begin{theorem}\label{t33}
Let $\{\lambda _j\}$ be the set of representatives of
$L^{\circ} /L$, then any irreducible $V_L^+$-module is
isomorphic to one of the following modules:$$V_L^{\pm},
V_{\lambda_j+L}(2\lambda_j \notin L),
V_{\lambda_j+L}^{\pm}(2\lambda_j \in L), (V_L^{T_{\chi}})^{\pm}.$$
\end{theorem}

Now we consider some decompositions of the modules for the vertex
operator algebras $M(1)^+$ and $V_L^+$. We denote $M_\h(1)$ for the
vertex operator algebra $M(1)$ associated with $\h$ and similarly
for the modules. It is clear that if $\h'$ is a subspace of $\h$
such that the restriction of the bilinear form on $\h$ to $\h'$ is
non-degenerate, then $M_{\h'}^+$ is a simple vertex operator
subalgebra of $M_{\h}^+$. Furthermore, if $\h=\h_1\bigoplus \h_2$
such that $(\h_1, \h_2)=0$, then the modules in Theorem \ref{t32}
viewed as $M_{\h_1}^+\otimes M_{\h_2}^+$-modules can be decomposed
as follows:
$$M_\h^+\cong (M_{\h_1}^+\otimes M_{\h_2}^+)\bigoplus
(M_{\h_1}^-\otimes M_{\h_2}^-),$$
$$M_\h^-\cong (M_{\h_1}^+\otimes M_{\h_2}^-)\bigoplus
(M_{\h_1}^-\otimes M_{\h_2}^+),$$
$$M_\h(1, \lambda)\cong M_{\h_1}(1, \lambda_1)\otimes
M_{\h_2}(1, \lambda_2),$$
$$M_\h(1)(\theta)^+\cong (M_{\h_1}(1)(\theta)^+\otimes M_{\h_2}(1)(\theta)^+)\bigoplus (M_{\h_1}(1)(\theta)^-\otimes M_{\h_2}(1)(\theta)^-),$$
$$M_\h(1)(\theta)^-\cong (M_{\h_1}(1)(\theta)^+\otimes M_{\h_2}(1)(\theta)^-)\bigoplus (M_{\h_1}(1)(\theta)^-\otimes M_{\h_2}(1)(\theta)^+).$$

Let $L=\mathbb{Z}\alpha$ be a positive definite even lattice of rank
one. Then all irreducible $V_L^+$-modules are decomposed into direct
sums of irreducible $M(1)^+$-modules as follows (cf. \cite{DG1} and
\cite{A1})
$$V_L^{\pm}\cong M(1)^\pm\bigoplus
_{m=1}^{\infty}{M(1,m\alpha)},$$
$$V_{\lambda+L}\cong \bigoplus _{m\in
\mathbb{Z}}{M(1, \lambda+m\alpha)},$$
$$V_{\frac{\alpha}{2}+L}^{\pm}\cong \bigoplus
_{m=0}^{\infty}{M(1, \frac{\alpha}{2}+m\alpha)},$$
$$(V_L^{T_i})^{\pm}\cong M(1)(\theta)^{\pm}, \ i=1,2$$
where $T_i$ is the one dimensional $\widehat{L}/K$-module $T_{\chi}$
with $\chi(e_{\al})=\pm 1$  respectively.

We also have the following result given in \cite{Y} and \cite{ABD}.
\begin{theorem}\label{n3.7}
 $V_L^+$ is
$C_2$-cofinite.
\end{theorem}
The following results were obtained in \cite{A2}.
\begin{theorem}\label{t3.3d}
The vertex operator algebra $V_L^+$ is rational, if $L$ is a
positive definite even lattice of rank one.
\end{theorem}
\begin{proposition}\label{t3.3c}
Let $L$ be a positive definite even lattice such that $V_{L}^{+}$ is
$C_2$-cofinite and $A(V_{L}^{+})$ is semisimple. Let $M^1, M^2$ be
irreducible $V_L^+$-modules. If the difference of the lowest weight
of $M^1$ and $M^2$ is not a nonzero integer, then
$$Ext_{V_L^+}^1(M^2, M^1) = 0.$$
\end{proposition}
\begin{remark}
In the next section, we will prove  that for a positive definite
even lattice $L$, the Zhu's algebra $A(V_{L}^{+})$ is semisimple. Thus
Proposition \ref{t3.3c} is true for any positive definite even lattice.
\end{remark}
 We now consider the
fusion rules for the vertex operator algebra $V_L^+$. For any
$\lambda \in L^{\circ} $ and a central character $\chi$ of
$\widehat{L}/K$ , let $\chi^{(\lambda)}$ be the central character
defined by $\chi^{(\lambda)}(a) = (-1)^{(\overline{a}, \lambda)}\chi
(a)$. We set $T_{\chi}^{(\lambda)} = T_{\chi^{(\lambda)}}$. We call
a triple $(\lambda,\mu,\nu)$ for $\lambda, \mu, \nu \in
L^{\circ} $ an admissible triple modulo $L$, if
$p\lambda+q\mu+r\nu \in L$ for some $p, q, r\in \{\pm 1\}$.

The following result on part of  fusion rules for the vertex
operator algebra $V_{L}^{+}$ when $L$ is of rank one comes from
\cite{A1}. This result will be used in Section 5 to deal with the rationality
of $V_L^+$ when $L$ has an orthogonal base.
\begin{theorem}\label{3.20}
Let $L=\Z\al$ be a positive definite even lattice and
$L^{\circ}/L=\{\la_0,\la_1,\la_2,\cdots, \la_k\}$ such that $\la_0=0$ and $\la_{k}=\al/2$.
Let $W^{i}, i=1,2,3$ be irreducible $V_{L}^{+}$-modules. Then

(1) \ the fusion rule of type $\left(\begin{tabular}{c}
$W^3$\\
$W^1$ $W^2$\\
\end{tabular}\right)$ is either 0 or 1.

(2) \ the fusion rule of type $\left(\begin{tabular}{c}
$W^3$\\
$W^1$ $W^2$\\
\end{tabular}\right)$  is non-zero if and only if $W^{i} (i=1,2,3)$
satisfy one of the following cases:

(i) \ $W^{1}=V_{L}^{+}$ and $W^{2}\cong W^{3}$.

(ii) \  $W^{1}=V_{L}^{-}$ and the pair $(W^{2},W^{3})$ is one of the
following pairs
$$(V_{L}^{\pm},V_{L}^{\mp}), \
(V_{\al/2+L}^{\pm},V_{\al/2+L}^{\mp}), \
(V_{L}^{T_{1},\pm},V_{L}^{T_{1},\mp}), \
(V_{L}^{T_{2},\pm},V_{L}^{T_{2},\mp}), \
(V_{\la_{i}+L},V_{\la_{i}+L}), \
$$
for $i=1,2,\cdots,k-1$.
\end{theorem}

The fusion rules for $V_L^+$ for any $L$ was obtained in \cite{ADL} and
will be exploited in Section 6. Recall the number $\pi_{\lambda,\mu}$ for $\lambda,\mu\in L^{\circ}$ from \cite{ADL}.
\begin{theorem}\label{t3.7}
Let $L$ be a positive definite even lattice , for any irreducible
$V_L^+$-modules $M^i$ $(i = 1, 2, 3),$ the fusion rule of type
$\left(\begin{tabular}{c}
$M^3$\\
$M^1$ $M^2$\\
\end{tabular}\right)$ is either 0 or 1. Furthermore, the fusion rule of type
$\left(\begin{tabular}{c}
$M^3$\\
$M^1$ $M^2$\\
\end{tabular}\right)$ with $M^1$ being one of $V_{\lambda+L}$ for
$(2\lambda\notin L)$ , $V_{\lambda+L}^+$ for $(2\lambda\in L)$ ,
$V_{\lambda+L}^-$ for $(2\lambda\in L)$ is 1 if and only if $M^i$ $(i=1, 2, 3)$ satisfy one of the following conditions:

(1) $M^1 = V_{\lambda+L}$ for $\lambda \in L^{\circ} $ such
that
$2\lambda \notin L$ and $M^2, M^3$is one of the following pairs:

$(V_{\mu+L}, V_{\nu+L})$ for $\mu, \nu \in L^{\circ} $ such
that $2\mu, 2\nu\notin L$ and $(\lambda, \mu, \nu)$ is an admissible
triple modulo $L.$

$(V_{\mu+L}^{\pm}, V_{\nu+L}), ((V_{\nu+L})', (V_{\mu+L}^{\pm})')$ for $\mu, \nu \in
L^{\circ} $ such that $2\mu\in L, 2\nu \notin L$ and
$(\lambda, \mu, \nu)$ is an admissible triple modulo $L.$

$(V_L^{T_{\chi}, \pm}, V_L^{T_{\chi}^{(\lambda)}, \pm}), (V_L^{T_{\chi}, \pm}, V_L^{T_{\chi}^{(\lambda)}, \mp})$
for any irreducible $\widehat{L}/K$-module $T_{\chi}$.

(2) $M^1 = V_{\lambda+L}^+$ for $\lambda \in L^{\circ} $
such that
$2\lambda \in L$ and $M^2, M^3$is one of the following pairs:\\
$(V_{\mu+L}, V_{\nu+L})$ for $\mu, \nu \in L^{\circ} $ such
that $2\mu, 2\nu\notin L$ and $(\lambda, \mu, \nu)$ is an admissible
triple modulo $L.$

$(V_{\mu+L}^{\pm},V_{\nu+L}^{\pm})$  for $\mu, \nu \in
L^{\circ} $ such that $2\mu , 2\nu \in L, \pi_{\lambda,
2\mu} = 1$  and $(\lambda, \mu, \nu)$ is an admissible triple modulo
$L.$

$(V_{\mu+L}^{\pm}, V_{\nu+L}^{\mp})$  for $\mu, \nu \in
L^{\circ} $ such that $2\mu , 2\nu \in L, \pi_{\lambda,
2\mu} = -1$  and $(\lambda, \mu, \nu)$ is an admissible triple
modulo $L.$

$(V_L^{T_{\chi}, \pm}, V_L^{T_{\chi}^{(\lambda)}, \pm}), ((V_L^{T_{\chi}^{(\lambda)}, \pm})', (V_L^{T_{\chi}, \pm})')$
for any irreducible $\widehat{L}/K$-module $T_{\chi}$ such that
$c_{\chi}(\lambda) = 1$.

$(V_L^{T_{\chi}, \pm}, V_L^{T_{\chi}^{(\lambda)},\mp}), ((V_L^{T_{\chi}^{(\lambda)}, \pm})', (V_L^{T_{\chi}, \mp})')$
for any irreducible $\widehat{L}/K$-module $T_{\chi}$ such that
$c_{\chi}(\lambda)=-1$.

(3) $M^1 = V_{\lambda+L}^-$ for $\lambda \in L^{\circ} $
such that
$2\lambda \in L$ and $M^2, M^3$is one of the following pairs:

$(V_{\mu+L}, V_{\nu+L})$ for $\mu, \nu \in L^{\circ} $ such
that $2\mu, 2\nu\notin L$ and $(\lambda, \mu, \nu)$ is an admissible
triple modulo $L.$

$(V_{\mu+L}^{\pm}, V_{\nu+L}^{\pm})$  for $\mu, \nu \in
L^{\circ} $ such that $2\mu , 2\nu \in L, \pi_{\lambda,
2\mu} = -1$  and $(\lambda, \mu, \nu)$ is an admissible triple
modulo $L.$

$(V_{\mu+L}^{\pm}, V_{\nu+L}^{\mp})$  for $\mu, \nu \in
L^{\circ} $ such that $2\mu , 2\nu \in L,\pi_{\lambda,2\mu}
= 1$  and $(\lambda, \mu, \nu)$ is an admissible triple modulo L
such that $c_{\chi}(\lambda)=1$.

$(V_L^{T_{\chi}, \pm}, V_L^{T_{\chi}^{(\lambda)}, \pm}), ((V_L^{T_{\chi}^{(\lambda)}, \pm})', (V_L^{T_{\chi}, \pm})')$
for any irreducible $\widehat{L}/K$-module $T_{\chi}$ such that
$c_{\chi}(\lambda) = -1$.

$(V_L^{T_{\chi}, \pm}, V_L^{T_{\chi}^{(\lambda)}, \mp}), ((V_L^{T_{\chi}^{(\lambda)}, \pm})', (V_L^{T_{\chi}, \mp})')$
for any irreducible $\widehat{L}/K$-module $T_{\chi}$ such that
$c_{\chi}(\lambda) = 1$.
\end {theorem}

Next we identify the contragredient modules of the irreducible
$V_L^+$-modules \cite{ADL}:
\begin{proposition}\label{n3.8}
The irreducible $V_L^+$-modules $V_L^{\pm}$ and $V_{\lambda + L}$
for $\lambda\in L^{\circ} $ with $2\lambda \notin L$ are
self dual. For $\lambda\in L^{\circ} $ with $2\lambda \in
L$, $V_{\lambda + L}^{\pm}$ are self dual if $2(\lambda, \lambda)$
is even, and $(V_{\lambda + L}^{\pm})'\cong V_{\lambda + L}^{\mp} $
if $2(\lambda, \lambda)$ is odd. Let $\chi$ be a central character
of $\widehat{L}/K$ such that $\chi(\kappa)=-1$, then the irreducible
modules $(V_L^{T_{\chi}, {\pm}})'$ are isomorphic to
$(V_L^{T_{\chi}', {\pm}})'$ respectively, where $\chi'$ is a central
character of $\widehat{L}/K$ defined by $\chi'(a) =
(-1)^{\frac{(\overline{a}, \overline{a})}{2}}\chi(a)$ for any $a\in
Z(\widehat{L}/K)$.
\end{proposition}
\section{Semisimplicity of $A(V_{L}^{+})$}
\def\theequation{4.\arabic{equation}}
\setcounter{equation}{0}

Motivated by Proposition \ref{t3.3c}, we prove the semisimplicity of
$A(V_L^+)$ for any positive definite even lattice $L$ in this
section. In the case that the rank of $L$ is 1, this result has
previously been obtained in \cite{DN2}. The semisimplicity of
$A(V_L^+)$ enables us to establish that if the two $V_L$-modules
have the same lowest weight then the extension of one module by the
other is always trivial.

First recall that  the irreducible $A(V_L^+)$-modules are the top
levels $W(0)$ of irreducible admissible $V_L^+$-modules $W$. So by
Theorem \ref{t33} (also see \cite{DN2}), we have
\begin{lem}
\label{l4.1c}
 The
irreducible $A(V_L^+)$-modules are given as follows:
\begin{eqnarray*}
&\charge{+}(0)=\C\1,\quad
\charge{-}(0)=\h(-1)\bigoplus(\bigoplus_{\alpha\in L_2}\C
(e^{\alpha}-e^{-\alpha})),\\
&\charlam{\lambda_i}(0)=\bigoplus_{\alpha\in\Delta(\lambda_i)}\C e^{\lambda_i+\alpha}\quad(2\lambda_i\notin L),\\
&\charlam{\lambda_i}^{\pm}(0)=\sum_{\alpha\in\Delta(\lambda_i)}\C(e^{\lambda_i+\alpha}\pm
\theta e^{\lambda_i+\alpha})\,
\quad(2\lambda_i\in L),\\
&\charge{T_{\chi},+}(0)=T_{\chi},\quad\charge{T_{\chi},-}(0)=\h(-1/2)\otimes
T_{\chi},
\end{eqnarray*}
 where $L_{2}=\{\al\in L| (\al, \al)=2\}$, $\h(-1)=\{h(-1)\1|h\in\h\}\subset M(1)$ and
$\h(-1/2)=\{h(-1/2)\1|h\in\h\}\subset M(1)(\theta).$
\end{lem}
Let $\{h_1,\cdots,h_d\}$ be an orthonormal basis of $\h.$ Recall
from \cite{DN2} and \cite{AD} the following vectors in $V_L^+$ for
$a,b=1,\cdots,d$ and $\alpha\in L$
\begin{align*}
S_{ab}&(m,n)=h_a(-m)h_b(-n),\\
E^u_{ab}&=5 S_{ab}(1,2)+25 S_{ab}(1,3)+36 S_{ab}(1,4)+16 S_{ab}(1,5)\,(a\neq b),\\
\bar{E}^u_{ba}&=S_{ab}(1,1)+14S_{ab}(1,2)+41S_{ab}(1,3)+44S_{ab}(1,4)+16 S_{ab}(1,5)\,(a\neq b),\\
E_{aa}^u&=E^{u}_{ab}E^{u}_{ba},\\
E^t_{ab}&=-16(3 S_{ab}(1,2)+14S_{ab}(1,3)+19S_{ab}(1,4)+8 S_{ab}(1,5))\,(a\neq b),\\
\bar{E}^t_{ba}&=-16(5S_{ab}(1,2)+18 S_{ab}(1,3)+21 S_{ab}(1,4)+8 S_{ab}(1,5))\,(a\neq b),\\
E_{aa}^t&=E^{t}_{ab}E^{t}_{ba},\\
\Lambda_{ab}&=45 S_{ab}(1,2)+190 S_{ab}(1,3)+240 S_{ab}(1,4)+96 S_{ab}(1,5),\\
E^{\a}&=e^{\a}+e^{-\a}.
\end{align*}
For $v\in V_L^+$, we denote $v+O(V_L^+)$ by $[v].$ Let $A^u$ and
$A^t$ be the linear subspace of $A(M(1)^+)$ spanned by
$E_{ab}^u+O(M(1)^+)$ and $E_{ab}^{t}+O(M(1)^+)$ respectively for
$1\leq a,\,b\leq d.$ Then $A^t$ and $A^u$ are two sided ideals of
$A(M(1)^+).$ Note that the natural algebra homomorphism from
$A(M(1)^+)$ to $A(V_L^+)$ gives embedding of $A^u$ and $A^t$ into
$A(V_L^+).$ We should remark that the $A^u$ and $A^t$ are
independent of the choice of the orthonormal basis
$\{h_1,\cdots,h_d\}.$

By Lemma 7.3 of \cite{AD} we know that
$$V_L^-(0)=\h(-1)\bigoplus(\sum_{\a\in L_2}\C[E^\a]\a(-1)),$$
where $L_2=\{\alpha\in L|(\a,\a)=2\}.$ Let
$L_{2}=\{\pm\al_{1},\cdots,\pm\al_{r},\pm\al_{r+1},\cdots,\pm\al_{r+l}\}$
be such that $\{\al_{1},\cdots,\al_{r}\}$ are linearly independent
and $\{\al_{r+1},\cdots,\al_{r+l}\} \ \subseteq
\bigoplus_{i=1}^{r}\Z_{+}\al_{i}.$ We can choose the orthonormal
basis $\{h_i | \ i=1,\cdots,d\}$ so that $h_{i}\in
\C\al_{1}+\cdots+\C\al_{i}$, for $i=1,\cdots,r.$ Then we have
$$\al_{i}(-1)=a_{i1}h_{1}(-1)+\cdots+a_{ii}h_{i}(-1), \
i=1,\cdots,r,$$
$$
\al_{j}(-1)=a_{j1}h_{1}(-1)+\cdots +a_{jr}h_{r}(-1), \
j=r+1,\cdots,r+l,$$ where $a_{ii}\neq 0, i=1,\cdots,r$. For
$i\in\{1,2,\cdots,l\}$, let $k_{i}$ be such that
$$a_{r+i,k_{i}}\neq 0, \ a_{r+i,k_{i}+1}=\cdots=a_{r+i,r}=0.$$

We know from \cite{AD} that $e^i=h_i(-1)$ for $i=1,\cdots,d$ and
$e^{d+j}=[E^{\a_j}]\alpha_j(-1)$ for $j=1,\cdots,r+l$ form a basis
of $V_L^-(0).$ We first construct a two-sided ideal of $A(V_L^+)$
isomorphic to $\End(V_L^-(0)).$ Recall $E^u_{ij}$ for
$i,j=1,\cdots,d.$ We now extend the definition of $E^u_{ij}$  to all
$i,j=1,\cdots,d+r+l$  and the linear span of $E^u_{ij}$ will be the
ideal of $A(V_L^+)$ isomorphic to $\End V_L^-(0)$ (with respect to
the basis $\{e^1,\cdots,e^{d+r+l}\}$ ).

For the notational convenience,  we also write $E^u_{i,j}$ for
$E^u_{ij}$ from now on. Define
$$
[E_{j,d+i}^{u}]=\frac{1}{4\epsilon(\al_{i},\al_{i})a_{ii}}[E_{ji}^{u}*E^{\al_{i}}],
\ i=1,\cdots,r,j=1,\cdots,d,$$
$$
[E_{j,d+r+i}^{u}]=\frac{1}{4\epsilon(\al_{r+i},\al_{r+i})a_{r+i,k_{i}}}[E_{j,k_{i}}^{u}*E^{\al_{r+i}}],
\ i=1,\cdots,l,j=1,\cdots,d.$$  Define
$$
[E_{d+i,j}^{u}]=\sum\limits_{k=1}^{r}a_{ik}[E^{\al_{i}}]*[E_{kj}^{u}],
\ i=1,\cdots,r+l,j=1,\cdots,d,$$ where  $a_{ij}=0$, for $1\leq
i<j\leq r$.
 Recall from \cite{DN2} and
\cite{AD} that $[E^u_{ab}]h_c(-1)=\delta_{c,b}h_a(-1)$ for
$a,b,c=1,\cdots,d.$
\begin{lem}\label{l5.1} The following holds:
$$[E^u_{ij}]e^k=\delta_{k,j}e^i,\ [E^u_{st}]e^k=\delta_{t,k}e^s$$
for $i,t=1,\cdots,d,$ $j,s=d+1,\cdots,d+r+l$ and $k=1,\cdots,d+r+l.$
\end{lem}

\pf Let $h\in {\mathfrak h}$ such that $(h,h)\ne 0.$ Then $\omega_h=
\frac{1}{2(h,h)}h(-1)^2$ is a Virasoro element with central charge
1. Note that $\omega_h\beta(-1)=\frac{(\beta,h)^4}{2(h,h)}h(-1)$ for
any $\beta\in {\mathfrak h}.$ For $\alpha\in L_{2}$ then
$[E^\alpha]*[E^\alpha]=4\e(\a,\a)[\omega_\alpha]$ in $A(V_L^+)$ by
Proposition 4.9 of \cite{AD}.  Then for $i=1,\cdots,d,
j=1,\cdots,r$, we have
\begin{eqnarray*}
& & [E_{i,d+j}^{u}]e^{d+j}=[E_{i,d+j}^{u}]([E^{\al_{j}}]\al_{j}(-1))\\
&
&\ \ \ =\frac{1}{4\epsilon(\al_{j},\al_{j})a_{jj}}([E_{ij}^{u}]*[E^{\al_{j}}]*[E^{\al_{j}}])\al_{j}(-1)\\
& &\ \ \ =\frac{1}{a_{jj}}[E_{ij}^{u}]\al_{j}(-1)=h_{i}(-1).
\end{eqnarray*}

Let $k\in\{1,\cdots,r+l\}$ such that $k\ne j.$ Then
\begin{eqnarray*}
& &[E_{i,d+j}^{u}]e^{d+k}=[E_{i,d+j}^{u}]([E^{\al_{k}}]\al_{k}(-1))\\
& &\ \ \
=\frac{1}{4\epsilon(\al_{j},\al_{j})a_{jj}}([E_{ij}^{u}]*[E^{\al_{j}}]*[E^{\al_{k}}])\al_{k}(-1).
\end{eqnarray*}
By Proposition 5.4 of \cite{AD}, we have
$$[E^{\al_{j}}]*[E^{\al_{k}}]=\sum\limits_{p}[v^{p}]*[E^{\al_{j}+\al_{k}}]*[w^{p}]+\sum\limits_{q}[x^{q}]*
[E^{\al_{j}-\al_{k}}]*[y^{q}],
$$
where $v^{p}, w^{p}, x^{q}, y^{q}\in M(1)^{+}$.
 Since $A^{u}$ is an ideal of
$A(M(1)^{+})$, we have $[E_{ji}^{u}]*[v^p], [E_{ji}^{u}]*[x^q]\in
A^{u}.$ By the proof of Proposition 7.2 of \cite{AD}, we know that
$A^{u}[E^{\al}]\al(-1)=0$, for any $\al\in L_{2}$. So by Lemma 7.1
and Proposition 7.2 of \cite{AD}, we have
$$[E_{i,d+j}^{u}]e^{d+k}=0, \ i=1,\cdots,d,
j=1,\cdots,r, k=1,\cdots,r+l, j\neq k.$$ It follows from the proof
of Proposition 7.2 of \cite{AD} that
$$E^u_{i,j+d}e^s=[E_{ij}^u]*[E^{\a_j}]h_s(-1)=0$$
for $s=1,\cdots,d$ as $[E^{\a_j}]h_s(-1)\in\sum_{p=1}^{r+l}\C
[e^{\a_{p}}-e^{-\a_{p}}].$ This completes the proof for
$E^u_{i,j+d}$ for $i=1,\cdots,d,$ $j=1,\cdots,r.$ The other cases
can be done similarly. \qed

Recall $H_a$ and $\omega_a=\omega_{h_a}$ for $a=1,\cdots,d$
 from \cite{DN2} and \cite{AD}. The following lemma collects some formulas
from Propositions 4.5, 4.6, 4.8 and 4.9 of \cite{AD}.
\begin{lem}\label{adl} For any indices $a,\,b,\,c,\,d$,
\begin{align}
&[\w_a]*[E^{u}_{bc}]=\delta_{ab}[E^{u}_{bc}],\label{e-1}\\
&[E^{u}_{bc}]*[\w_a]=\delta_{ac}[E^{u}_{bc}].\label{e-2}\\
&[E_{ab}^{u}]*[E_{cd}^{t}]=[E_{cd}^{t}]*[E_{ab}^{u}]=0,
\label{equ1}\\
&[\La_{ab}]*[E_{cd}^{u}]=[\La_{ab}]*[E_{cd}^{t}]=[E_{cd}^{u}]*[\La_{ab}]=[E_{cd}^{t}]*[\La_{ab}]=0
\ (a\neq b).\label{equ2}
\end{align}
For distinct $a,b,c,$
\begin{align}
&\left(70 [H_a]+1188[\w_a]^2-585 [\w_a]+27\right)*[H_a]=0,\label{equation1}\\
&([\w_a]-1)*\left([\w_a]-\frac{1}{16}\right)*\left([\w_a]-\frac{9}{16}\right)*[H_a]=0,\label{equation2}\\
&-\frac{2}{9}[H_a]+\frac{2}{9}[H_b]=2[E_{aa}^u]-2[E_{bb}^u]+\frac{1}{4}[E_{aa}^t]-\frac{1}{4}[E_{bb}^t],
\label{equation3}\\
\begin{split}\label{equation4}
&-\frac{4}{135}(2[\w_a]+13)*[H_a]+\frac{4}{135}(2[\w_b]+13)*[H_b]\\
&\ \ \
=4([E_{aa}^u]-[E_{bb}^u])+\frac{15}{32}([E_{aa}^t]-[E_{bb}^t]),
\end{split}\\
&[\w_b]*[H_a]=-\frac{2}{15}([\w_a]-1)*[H_a]+\frac{1}{15}([\w_b]-1)*[H_b],\label{equation5}\\
&[\Lambda_{ab}]^2=4[\w_a]*[\w_b]-\frac{1}{9}([H_a]+[H_b])-([E_{aa}^u]+[E_{bb}^u])-\frac{1}{4}
([E_{aa}^t]+[E_{bb}^t]),\label{equation6}\\
&[\Lambda_{ab}]*[\Lambda_{bc}]=2[\w_b]*[\Lambda_{ac}].\label{equation7}
\end{align}
For $\alpha\in L$ such that $(\alpha,\alpha)=2k\ne 2,$
\begin{align}
&[H_{\alpha}]*[E^{\alpha}]=\frac{18(8k-3)}{(4k-1)(4k-9)}\left([\w_{\alpha}]-\frac{k}{4}\right)\left([\w_{\alpha}]-
\frac{3(k-1)}{4(8k-3)}\right)[E^{\alpha}],\label{equation8}\\
&\left([\w_{\alpha}]-\frac{k}{4}\right)\left([\w_{\alpha}]-\frac{1}{16}\right)\left([\w_{\alpha}]-
\frac{9}{16}\right)[E^{\alpha}]=0.\label{equation9}
\end{align}
If $\alpha\in L_2,$
\begin{align}
&[E^\alpha]*[E^\alpha]=4\e(\a,\a)[\w_\alpha],\label{equation10}\\
&[H_\alpha]*[E^\alpha]+[E^\alpha]*[H_\alpha]=-12[\w_\alpha]*\left([\w_{\alpha}]-\frac{1}{4}\right)*[E^{\alpha}],
\label{equation11}\\
&([\w_{\alpha}]-1)*\left([\w_{\alpha}]-\frac{1}{4}\right)*\left([\w_{\alpha}]-\frac{1}{16}\right)*\left([\w_{\alpha}]-
\frac{9}{16}\right)*[E^{\alpha}]=0.\label{equation12}
\end{align}
For any $\al\in L$,
\begin{align}
I^{t}*[E^{\al}]=[E^{\al}]*I^{t},\label{equation13}
\end{align}
where $I^{t}$ is the identity of the simple algebra $A^{t}$.
\end{lem}

\begin{lem}\label{l5.2} For any $\al\in L_{2}$, we have
$$A^u*[E^{\al}]*A^u=0.
$$
\end{lem}

\pf Let $\al\in L_{2}$ and $\{h_{1},\cdots,h_{d}\}$ be  an
orthonormal basis of ${\f h}$ such that $h_{1}\in\C\al$. ($A^u$ is
independent of the choice of orthonormal basis.) By
(\ref{e-1})-(\ref{e-2}) and (\ref{equation12}), we have
$[E^{\alpha}]=f([\omega_{\a}])*[E^{\alpha}]=[E^{\alpha}]*f([\omega_{\a}])$
for some polynomial $f(x)$ with $f(0)=0.$ Note that
$\omega_{\alpha}=\omega_1.$ By (\ref{e-1})-(\ref{e-2}), we  only
need to prove that
$$
[E_{i1}^{u}]*[E^{\al}]*[E_{1s}^{u}]=0, \ i,s=1,2,\cdots,d.$$ Let
$a=1$, $b\neq 1$ in (\ref{equation5}). Multiplying (\ref{equation5})
by $[E_{i1}^{u}]$ on left and using (\ref{e-2}) and (\ref{equ1}), we
have
$$[E_{i1}^{u}]*[H_{b}]=0, \ b\neq 1.$$
Then setting $a=1$, $b\neq 1$ in (\ref{equation3}) and multiplying
(\ref{equation3}) by $[E_{i1}^{u}]$ on left yields
$$[E_{i1}^{u}]*[H_{1}]=-9[E_{i1}^{u}].$$
Let $a=1$, $b\neq 1$ in (\ref{equation6}). Multiplying
(\ref{equation6}) by $[E_{1s}^{u}]$ on right and using (\ref{e-1})
and (\ref{equ2}), we have
$$
-\frac{1}{9}[H_{1}]*[E_{1s}^{u}]-\frac{1}{9}[H_{b}]*[E_{1s}^{u}]=[E_{1s}^{u}].$$
On the other hand, multiplying (\ref{equation3}) by $[E_{1s}^{u}]$
on right yields
$$
-\frac{1}{9}[H_{1}]*[E_{1s}^{u}]+\frac{1}{9}[H_{b}]*[E_{1s}^{u}]=[E_{1s}^{u}].$$
Comparing the above two formulas, we have
$$[H_{1}]*[E_{1s}^{u}]=-9[E_{1s}^{u}], \
[H_{a}]*[E_{1s}^{u}]=0, \ a\neq 1.$$ So
$$
[E_{i1}^{u}]*[H_{1}]*[E^{\al}]*[E_{1s}^{u}]+[E_{i1}^{u}]*[E^{\al}]*[H_{1}]*
[E_{1s}^{u}]=-18[E_{i1}^{u}]*[E^{\al}]*[E_{1s}^{u}].$$ But by
(\ref{e-1})-(\ref{e-2}) and (\ref{equation11})
 we have
$$
[E_{i1}^{u}]*[H_{1}]*[E^{\al}]*[E_{1s}^{u}]+[E_{i1}^{u}]*[E^{\al}]*[H_{1}]*
[E_{1s}^{u}]=-9[E_{i1}^{u}]*[E^{\al}]*[E_{1s}^{u}].$$ This implies
that $ [E_{i1}^{u}]*[E^{\al}]*[E_{1s}^{u}]=0,$ as required.
 \qed

We now define $E^u_{i,j}$ for all $i,j=1,\cdots,d+r+l.$  Set
$$
[E_{d+i,d+j}^{u}]=[E_{d+i,1}^{u}]*[E_{1,d+j}^{u}], \
i,j=1,\cdots,r+l.$$ It is easy to see that
$[E_{d+i,1}^{u}]*[E_{1,d+j}^{u}]=[E_{d+i,k}^{u}]*[E_{k,d+j}^{u}]$,
$k=2,\cdots,d$.

Denote by $A_{L}^u$ the subalgebra of $A(V_{L}^{+})$ generated by
$$\{[E_{ij}^{u}],[E_{d+p,j}^{u}],[E_{i,d+p}^{u}]|i,j=1,\cdots,d, \
p=1,\cdots,r+l\}.$$
 From Lemma \ref{l5.2}, (\ref{equation10}) and the
definition of $[E_{ij}^{u}]$, $i,j=1,\cdots,d+r+l$, we can easily
deduce the following result.

\begin{lem}\label{l5.4} $A_{L}^u$ is a matrix algebra over $\C$ with basis  $\{[E_{ij}^{u}]|i,j=1,\cdots,d+r+l\}$
such that
$$
[E_{ij}^{u}]*[E_{ks}^{u}]=\de_{j,k}[E_{is}^{u}], \
[E_{ij}^{u}]e^{k}=\de_{j,k}e^{i}, \ i,j,k,s=1,2,\cdots,d+r+l.$$
\end{lem}

\begin{lem}\label{l5.3} Let $\al\in L$ be such that $\al\notin L_{2}$, then
$$[E^{\al}]*A^u=0.$$
\end{lem}

\pf Let $\{h_{1},\cdots,h_{d}\}$ be  an orthonormal basis of ${\f
h}$ such that $h_{1}\in\C\al$. If $|\al|^{2}=2k$ and $k\neq 4,$ the
lemma follows from (\ref{e-1})-(\ref{e-2}) and (\ref{equation9}).

If $|\al|^{2}=8,$ by (\ref{e-1})-(\ref{e-2}) and (\ref{equation9})
we have
$$
[E_{ab}^{u}]*[E^{\al}]=[E^{\al}]*[E_{ba}^{u}]=0, $$ for all $1\leq
a,b\leq d$ and $b\neq 1$. By (\ref{e-1}) and (\ref{equation8}), we
have
\begin{equation}\label{ea5.20}
[E_{a1}^{u}]*[H_{1}]*[E^{\al}]=0.
\end{equation}
On the other hand, for $a\neq 1$, by
(\ref{equation3})-(\ref{equation4}) and (\ref{equ1}), we have
$$
-\frac{2}{9}[E_{a1}^{u}]*([H_{a}]-[H_{1}])*[E^{\al}]=-2[E_{a1}^{u}]*[E^{\al}],$$
$$
\frac{4}{135}[E_{a1}^{u}]*(-13[H_{a}]+15[H_{1}])*[E^{\al}]=-4[E_{a1}^{u}]*[E^{\al}].$$
Therefore by (\ref{ea5.20}), we have
$$\frac{1}{9}[E_{a1}^{u}]*[H_{a}]*[E^{\al}]=[E_{a1}^{u}]*[E^{\al}], \ a\neq 1,$$
$$\frac{13}{135}[E_{a1}^{u}]*[H_{a}]*[E^{\al}]=[E_{a1}^{u}]*[E^{\al}], \ a\neq 1.$$
This means that
$$
[E_{a1}^{u}]*[E^{\al}]=0.$$ Since
$[H_{\al}]*[E^{\al}]=[E^{\al}]*[H_{\al}]$, we similarly have
$$[E^{\al}]*[E_{1a}^{u}]=0.$$
This completes the proof. \qed

\begin{lem}\label{l5.5} $A_{L}^{u}$ is an ideal of $A(V_{L}^{+})$.
\end{lem}

\pf By Proposition 5.4  of \cite{AD}, (\ref{equ1}),
(\ref{equation13}) and Lemmas \ref{l5.2}-\ref{l5.3}, it is enough to
prove that $[E^{\al_{i}}]*[E_{jk}^{u}],
[E_{jk}^{u}]*[E^{\al_{i}}]\in A_{L}^u$,
  $j,k=1,\cdots,d, \
 i=1,\cdots,r+l$.

 Let $\al\in L_{2}$. For convenience, let $a_{ij}=0$, for
$1\leq i<j\leq r$ and $k_{i}=i$ for $1\leq i\leq r$. Since
$\al_{i}(-1)=\sum\limits_{k=1}^{d}a_{ik}h_{k}(-1)$, we have
$$\omega_{\al_{i}}=\frac{1}{4}\sum\limits_{k=1}^{d}a_{ik}^{2}h_{k}(-1)^{2}+\frac{1}{4}\sum\limits_{p\neq
q}a_{ip}a_{iq}h_{p}(-1)h_{q}(-1).$$ Recall from \cite{AD} that
$$
[S_{ab}(1,1)]=[E_{ab}^{u}]+[E_{ba}^{u}]+[\La_{ab}]+\frac{1}{2}[E_{ab}^{t}]+\frac{1}{2}[E_{ba}^{t}],
\ a\neq b.$$ So from (\ref{equ1}) and (\ref{equ2}) we have
$$
[E_{jk}^{u}]*[\omega_{\al_{i}}]=\frac{1}{2}a_{ik}^{2}[E_{jk}^{u}]+\frac{1}{2}\sum\limits_{p\neq
k}a_{ik}a_{ip}[E_{jp}^{u}], \ j,k=1,\cdots,d, \ i=1,\cdots,r+l.$$
Then it can  easily be deduced that
$$(a_{ik}[E_{jk_{i}}^{u}]-a_{ik_{i}}[E_{jk}^{u}])*[\omega_{\al_{i}}]=0,
\ j,k=1,\cdots,d, \ i=1,\cdots,r+l.$$ Then by (\ref{equation12}), we
have
$$(a_{ik}[E_{jk_{i}}^{u}]-a_{ik_{i}}[E_{jk}^{u}])*[E^{\al_{i}}]=0,
\ j,k=1,\cdots,d, \ i=1,\cdots,r+l.$$ So
$[E_{jk}^{u}]*[E^{\al_{i}}]\in A_{L}^{u}$, for all $j,k=1,\cdots,d,
i=1,\cdots,r+l.$ Similarly, we have
$$[E^{\al_{i}}]*[E_{jk}^{u}]=0, \ k=1,\cdots,d, \ i=1,\cdots,r+l, \ j=r+1,\cdots,d,$$
$$
[E^{\al_{i}}]*[\sum\limits_{b=1}^{r}a_{kb}E_{bj}^{u}]=\frac{(\al_{i},\al_{k})}{2}[E^{\al_{i}}]*[\sum\limits_{b=1}^{r}
a_{ib}E_{bj}^{u}], $$ for  $j=1,\cdots,d, \ k=1,\cdots,r, \
i=1,\cdots,r+l.$ Since
 both $\{\al_{1},\cdots,\al_{r}\}$ and $\{h_{1},\cdots,h_{r}\}$
 are linearly independent, it follows that for each
 $i=1,\cdots,r, \ j=1,\cdots,d$,
 $[E_{ij}^{u}]$ is a linear combination of $a_{11}[E_{1j}^{u}]$,
 $[a_{21}E_{1j}+a_{22}E_{2j}^{u}], \cdots,
 [a_{r1}E_{1j}^{u}+\cdots+a_{rr}E_{rj}^{u}]$. Therefore
 $[E^{\al_{i}}]*[E_{jk}^{u}]\in A_{L}^u$, $j,k=1,\cdots,d, \
 i=1,\cdots,r+l$. \qed

\vskip 0.2cm
 For $0\neq \al\in L$, let $\{h_{1},\cdots,h_{d}\}$ be  an
orthonormal basis of ${\f h}$ such that $h_{1}\in\C\al$. Define
$$[B_{\al}]=2^{|\al|^{2}-1}([I^{t}]*[E^{\al}]-\frac{2|\al|^{2}}{2|\al|^{2}-1}[E_{11}^{t}]*[E^{\al}]),$$
and $[B_{0}]=[I^{t}]$ (see formula (6.5) of \cite{AD}).
\begin{lem}\label{l5.6}
For $\al\in L$, $[E_{ij}^{t}]\in A^{t}$,
$[B_{\al}]*[E_{ij}^{t}]=[E_{ij}^{t}]*[B_{\al}]$.
\end{lem}

\pf It is enough to prove that
$$
[B_{\al}]*[E_{ij}^{t}]=[E_{ij}^{t}]*[B_{\al}],$$ for $i=1$ or $j=1$.
By the definition of $[E_{ab}^{t}]$ and the fact that
$$[I^{t}]*[E^{\al}]=[E^{\al}]*[I^{t}]$$
and
$$
[I^{t}]*[\La_{ab}]=[I^{t}]*[E_{ab}^{u}]=0, \ a\neq b,$$ we have
$$
[B_{\al}]*[E_{ab}^{t}]=[B_{\al}]*(-[S_{ab}(1,1)]-2[S_{ab}(1,2)]),$$
$$
[E_{ab}^{t}]*[B_{\al}]=(-[S_{ab}(1,1)]-2[S_{ab}(1,2)])*[B_{\al}].$$
 Let $b\neq 1$. Similar to
the proof of Lemma 7.5 of \cite{AD}, we have
\begin{eqnarray*}
& & \ \ \
(2|\al|^{2}-1)([E_{1b}^{t}]+3[E_{1b}^{u}]+[\La_{1b}])*[E^{\al}]+[E^{\al}]*([E_{1b}^{t}]+3[E_{1b}^{u}]+[\La_{1b}])\\
& &
=-([E_{b1}^{t}]-[E_{1b}^{u}]+[\La_{1b}])*[E^{\al}]-(2|\al|^{2}-1)[E^{\al}]*([E_{b1}^{t}]-[E_{1b}^{u}]+[\La_{1b}])
\end{eqnarray*}
\begin{eqnarray*}
&& \ \ \ (2|\al|^{2}-1)(\frac{1}{16}[E_{1b}^{t}]*[E^{\al}]+[\omega_{b}]*[\La_{1b}]*[E^{\al}])\\
&&\ \ \ \ \ \ \ \ +\frac{1}{16}[E^{\al}]*[E_{1b}^{t}]+[E^{\al}]*[\omega_{b}]*[\La_{1b}]\\
&&
\ \ \ =\frac{9}{16}[E_{b1}^{t}]*[E^{\al}]\!+\![\omega_{b}]*[\La_{1b}]*[E^{\al}]\\
&&\ \ \ \ \ \ \ \   -(2|\al|^{2}-1)(\frac{9}{16}[E^{\al}]*[E_{b1}^{t}]+[E^{\al}]*[\omega_{b}]*[\La_{1b}]).
\end{eqnarray*}
So we have
$$(2|\alpha|^{2}-1)[E_{1b}^{t}]*[E^{\al}]=-[E^{\al}]*[E_{1b}^{t}]+x,$$
where $x\in
A_{L}^{u}+\C[\La_{1b}]*[E^{\al}]+\C[E^{\al}]*[\La_{1b}]+\C[\omega_{b}]*[\La_{1b}]*[E^{\al}]
+\C[E^{\al}]*[\La_{1b}]*[\omega_{b}]$. Since $y*x=0$ for any $y\in
A^{t}$, we have
\begin{eqnarray*}
& &  \ \ \
[B_{\al}]*[E_{1b}^{t}]\\
 &
 &=2^{|\al|^{2}-1}\left(-(2|\al|^{2}-1)[E_{1b}^{t}]*[E^{\al}]+2|\al|^{2}[E_{11}^{t}]*[E_{1b}^{t}]*[E^{\al}]\right)\\
 & &=2^{|\al|^{2}-1}[E_{1b}^{t}]*[E^{\al}]=[E_{1b}^{t}]*[B_{\al}].
\end{eqnarray*}
Similarly,
$$[B_{\al}]*[E_{b1}^{t}]=[E_{b1}^{t}]*[B_{\al}],$$
completing the proof. \qed

\begin{lem}\label{l5.7} $A_{L}^{t}$ is an ideal of $A(V_{L}^{+})$ and
$A_{L}^{t}\cong A^{t}\otimes_{\C} {\C}[\widehat{L}/ K]/J$, where
${\C}[\widehat{L}/K]$ is the group algebra of $\widehat{L}/ K$ and
$J$ is the ideal of ${\C}[\widehat{L}/ K]$ generated by $\kappa
K+1.$
\end{lem}

\pf By Proposition 5.4 of \cite{AD} and Lemmas
\ref{l5.5}-\ref{l5.6}, it is easy to check that $A_{L}^{t}$ is an
ideal of $A(V_{L}^{+})$. Similar to the proof of Proposition 7.6 of
\cite{AD}, we have
$$
[B_{\al}]*[B_{\be}]=\epsilon(\al,\be)[B_{\al+\be}],$$ for
$\al,\be\in L$ where $\epsilon(\al,\be)$ is understood to be $\pm 1$
by identifying $\kappa$ with $-1.$ Then the lemma follows from
Proposition 7.6 of \cite{AD} and Lemma \ref{l5.6}. \qed

\vskip 0.2cm It is clear that $A_{L}^{u}\cap A_{L}^{t}=0$. Let
$$
\bar{A}(V_{L}^{+})=A(V_{L}^{+})/(A_{L}^{u}\oplus A_{L}^{t}),
$$
and for $x\in {A}(V_{L}^{+})$, we still denote the image of $x$ in
$\bar{A}(V_{L}^{+})$ by $x$.

\begin{lem}\label{l5.8}
In $\bar{A}(V_{L}^{+})$, we have
\begin{equation}\label{ea5.22}
[H_{a}]=[H_{b}], \ 1\leq a,b\leq d,
\end{equation}
\begin{equation}\label{ea5.23}
([\omega_{a}]-\frac{1}{16})*[H_{a}]=0, \ 1\leq a\leq d,
\end{equation}
\begin{equation}\label{ea5.24}
\frac{128}{9}[H_{a}]*\frac{128}{9}[H_{a}]=\frac{128}{9}[H_{a}], \
1\leq a\leq d,
\end{equation}
\begin{equation}\label{ea5.25}
[\La_{ab}]*[H_{c}]=0, \ 1\leq a,b,c\leq d, \ a\neq b.
\end{equation}
\end{lem}

\pf (\ref{ea5.22}) follows from (\ref{equation3})  and
(\ref{ea5.23}) follows from (\ref{equation4}) and (\ref{equation5}).
Then from (\ref{equation1}) we can get (\ref{ea5.24}). By
(\ref{equation6}), we have
$$
[\La_{ab}]^{2}*[H_{c}]=0, \ a\neq b.
$$
If $d\geq 3$, then by (\ref{ea5.22}) we can let $c\neq a$, $c\neq
b$. So by (\ref{equation7}) and (\ref{ea5.23}),
\begin{eqnarray*}
& & \   \ \ [\La_{ab}]*[H_{c}]=16[\La_{ab}]*[\omega_{c}]*[H_{c}]\\
& &=8[\La_{ac}]*[\La_{cb}]*[H_{c}]=128[\La_{ac}]*[\La_{cb}]*[\omega_{a}]*[H_{a}]\\
& &=64[\La_{ac}]*[\La_{ca}]*[\La_{ab}]*[H_{a}]=0.
\end{eqnarray*}
If $d=2$. Notice that $[\La_{ab}]=[S_{ab}(1,1)]$. By Remark 4.1.1 of
\cite{DN2} and the fact that
$[\omega_{a}*S_{ab}(m,n)]=[S_{ab}(m,n)*\omega_{a}]$ in
$\bar{A}(V_{L}^{+})$ for $m,n\geq 1$, we have
\begin{equation}[S_{ab}(m+1,n)]+[S_{ab}(m,n)]=0.\label{ea30}
\end{equation} By the proof of Lemma 6.1.2 of \cite{DN2}, we know that
\begin{equation}
[H_{a}]=-9[S_{aa}(1,3)]-\frac{17}{2}[S_{aa}(1,2)]+\frac{1}{2}[S_{aa}(1,1)].\label{ea31}
\end{equation}
Direct calculation yields
$$
[S_{ab}(1,1)]*[S_{aa}(1,3)]=h_{b}(-1)h_{a}(-3)h_{a}(-1)^{2},$$
$$
[S_{ab}(1,1)]*[S_{aa}(1,2)]=h_{b}(-1)h_{a}(-2)h_{a}(-1)^{2},$$
$$
[S_{ab}(1,1)]*[S_{aa}(1,1)]=h_{b}(-1)h_{a}(-1)h_{a}(-1)^{2}.$$ Here
we have used (\ref{ea30}). Then (\ref{ea5.25}) immediately follows
from Lemma 4.2.1 of \cite{DN2}, (\ref{ea30}) and (\ref{ea31}).
 The proof is
complete. \qed

\vskip 0.2cm
 For $0\neq \al\in L$, let $\{h_{1},\cdots,h_{d}\}$ be
an orthonormal basis of ${\f h}$ such that $h_{1}\in\C\al$. Define
$$
[\bar{B}_{\al}]=2^{|\al|^{2}-1}\frac{128}{9}[H_{1}]*[E^{\al}]. $$ We
also set $[\bar{B}_{0}]=\frac{128}{9}[H_{1}]$.

\begin{lem}\label{l5.9}
 The subalgebra $A_{H}$ of $\bar{A}(V_{L}^{+})$ spanned by
$[\bar{B}_{\al}]$, $\al\in L$ is an ideal of $\bar{A}(V_{L}^{+})$
isomorphic to ${\C}[\widehat{L}/K]/J$.
\end{lem}

Let
$$
\widehat{A}(V_{L}^{+})=\bar{A}(V_{L}^{+})/A_{H}.$$

\begin{lem}\label{laa}
 Any $\widehat{A}(V_{L}^{+})$-module is
completely reducible. That is, $\widehat{A}(V_{L}^{+})$ is a
semisimple associative algebra.
\end{lem}

\pf Let $M$ be an $\widehat{A}(V_{L}^{+})$-module. For $\al\in L$,
by \cite{DN2} $M$ is a direct sum of irreducible
$A(V_{\Z\al}^{+})$-modules. Following the proof of Lemma 6.1 of
\cite{AD} one can prove that the image of any vector from $M(1)^{+}$
in $\widehat{A}(V_{L}^{+})$ is semisimple on $M$. By Table 1 of
\cite{AD}, we can assume that
$$
M=\bigoplus_{\la\in{\f h/(\pm 1)}}M_{\la},$$ where $M_{\la}=\{w\in
M| \ [\frac{1}{2}h(-1)^{2}1]w=\frac{1}{2}(\la,h)^{2}w, h\in{\f
h}\}$. So $\omega_{a}w=\frac{1}{2}(\la,h_{a})^{2}w$, for $w\in
M_{\la}$. By (\ref{equation6}) and (\ref{equation7}), we have
$$\La_{ab}w=(\la,h_{a})(\la,h_{b})w, $$
for $a\neq b$, $w\in M_{\la}.$ For any $u\in M_{\la}$, $\la\neq 0$,
set $M(u)=\sum\limits_{\al\in L}\C[E^{\al}]u$. By
(\ref{equation8})-(\ref{equation9}) and
(\ref{equation11})-(\ref{equation13}), if $[E^{\al}]u\neq 0$, then
$\al\in\De(\la)$ or $-\al\in \De(\la)$, where $\De(\la)=\{\al\in L|
\ |\la+\al|^{2}=|\la|^{2}\}$. So
$$M(u)=\bigoplus_{\al\in\De(\la)}\C[E^{\al}]u.$$
Since $L$ is positive-definite, there are finitely many $\al\in L$
which belong to $\De(\la)$. Thus $M(u)$ is finite-dimensional.
Similar to the proof of  Lemma 6.4 of \cite{AD}, we can deduce that
$\La_{ab}M(u)\subseteq M(u)$, $\omega_{a}M(u)\subseteq M(u)$. By
Proposition 5.4 of \cite{AD},
$[E^{\al}]*[E^{\be}]=[x]*[E^{\al+\be}]$ for some $x\in M(1)^{+}$. We
deduce that $M(u)$ is an $\widehat{A}(V_{L}^{+})$-submodule of $M$.
Suppose $[E^{\al}]u\neq 0$, for some $\al\in\De(\la)$. If
$(\al,\al)=2$, then by (\ref{equation10}), we have $0\neq
[E^{\al}][E^{\al}]u\in{\C}u$. If $(\al,\al)=2k\neq 2$. Let
$\{h_{1},\cdots,h_{d}\}$ be an orthonormal basis of ${\f h}$ such
that $h_{1}\in\C\al$. By the fact that
$[H_{1}]=[J_{1}]+[\omega_{1}]-4[\omega_{1}^{2}]=0$ and
(\ref{equation9}) we know that $[\omega_{1}]u=\frac{k}{4}u.$ Then by
Lemma 5.5 of \cite{DN2}, we have
$$
[E^{\al}][E^{\al}]u=\frac{2k^{2}}{(2k)!}(k^{2}-1)(k^{2}-2^{2})\cdots(k^{2}-(k-1)^{2})u\neq
0.$$
  Therefore $M(u)$ is
irreducible. We prove that $M$ is a direct sum of finite-dimensional
irreducible $\widehat{A}(V_{L}^{+})$-modules. \qed

\begin{theorem}\label{t5.10} $A(V_L^+)$ is a finite dimensional semisimple associative
algebra.
\end{theorem}
\pf Clearly $A(V_L^+)$ is finite dimensional. By Lemmas \ref{l5.4}, \ref{l5.5}, \ref{l5.7} we know that $A_L^u\oplus A_L^t$ is a semisimple ideal of $A(V_L^+).$
Thus $A(V_L^+)$ is semisimple if and only if $\bar A(V_L^+)=A(V_L^+)/(A_L^u\oplus A_L^t)$ is semisimple. By Lemma \ref{l5.9}, $A_H$ is a semisimple ideal of
$\bar A(V_L^+).$ So $\bar A(V_L^+)$ is semisimple if and only if $\widehat{A}(V_{L}^{+})=\bar A(V_L^+)/A_H$ is semisimple. The result now follows
from Lemma \ref{laa}. \qed

\section{Rationality of $V_L^+$  when $L$ has an orthogonal base}

In this section we assume that $L$ has an orthogonal base
$\{\beta_i, 1\leq i \leq d\}$ in the sense that $(\beta_i, \beta_j) =
0$, for $i\neq j$. Then we have $L = \bigoplus_{i =
1}^d{\mathbb{Z}\beta_i}$ and this induces the relations
$$\otimes_{i = 1}^{d}{V_{\mathbb{Z}\beta_i}^+}\subseteq
(V_{\bigoplus_{i = 1}^d{\mathbb{Z}\beta_i}})^+ =  V_L^+$$ between
vertex operator algebras.  By Theorems \ref{t3.3d} and \ref{t2.6},
$\otimes_{i = 1}^{d}{V_{\mathbb{Z}\beta_i}^+}$ is a rational vertex operator subalgebra
of $V_L^+.$ This is a crucial fact in our discussion of rationality of $V_L^+$ in this section.

The following lemma is trivial:
\newtheorem{name5}{Lemma}[section]
\begin{lem}\label{lc5.1}
The vertex operator algebras $ \otimes_{i=
1}^{d}{V_{\mathbb{Z}\beta_i}}$ and $ V_{\bigoplus_{i=
1}^{d}{\mathbb{Z}\beta_i}}$ are isomorphic. \end{lem}

 By Lemma \ref{lc5.1}, a $
V_{\bigoplus_{i= 1}^{d}{\mathbb{Z}\beta_i}}$-module $
V_{\bigoplus_{i= 1}^{d}{\mathbb{Z}\beta_i}+c_i{\beta_i}}$ can be
viewed as a $ \otimes_{i= 1}^{d}{V_{\mathbb{Z}\beta_i}}$-module, and
is isomorphic to $ \otimes_{i=
1}^{d}{V_{\mathbb{Z}\beta_i+c_i\beta_i}}$.  On the other hand, by
Theorem \ref{t2.6} and Theorem \ref{t3.3d}, $\otimes_{i=
1}^{d}{V_{\mathbb{Z}\beta_i}^+}$ is rational, then any admissible
$V_L^{+}$-module is completely reducible as an admissible
$\otimes_{i= 1}^{d}{V_{\mathbb{Z}\beta_i}^+}$-module. Thus we may
decompose all the irreducible $V_L^{+}$-modules as a direct sum of
irreducible $\otimes_{i= 1}^{d}V_{\mathbb
{Z}\beta_{i}}^{+}$-modules.
\begin{lem}\label{l5.2c} The following are $\otimes_{i= 1}^{d}V_{\mathbb
{Z}\beta_{i}}^{+}$-module isomorphisms.

(1) \ $V_L^+$ is isomorphic to a direct sum of
$$V_{\mathbb{Z}\beta_1}^{\epsilon_1}\otimes\cdots\otimes
V_{\mathbb{Z}\beta_d}^{\epsilon_d}$$ with sign ${\epsilon_i=
\{\pm\}, 1\leq i\leq d}$ such that the number of $i$ with
$\epsilon_i= -$ is even.

(2) $V_L^-$ is isomorphic to a direct sum of
$$V_{\mathbb{Z}\beta_1}^{\epsilon_1}\otimes\cdots\otimes
V_{\mathbb{Z}\beta_d}^{\epsilon_d}$$ with sign ${\epsilon_i=
\{\pm\}, 1\leq i\leq d}$ such that the number of $i$ with
$\epsilon_i= -$ is odd.

(3) $V_{\lambda_j+L}$ for $ \lambda_j=
k_1\beta_1+\cdots+k_d\beta_d\in L^{\circ} /L, k_1, \cdots,
k_d\in \mathbb{C}$ and $2\lambda_j\notin L$ is isomorphic to
$$V_{\mathbb{Z}\beta_1+k_1\beta_1} \otimes\cdots\otimes
V_{\mathbb{Z}\beta_d+k_d\beta_d}, $$ where
$V_{\mathbb{Z}\beta_i+k_i\beta_i}=V_{\mathbb{Z}\beta_i+k_i\beta_i}^{+}\bigoplus
V_{\mathbb{Z}\beta_i+k_i\beta_i}^{-}$ for $i$ such that
$2k_{i}\be_{i}\in {\mathbb Z}\be_{i}$.

(4) $V_{\lambda_j+L}^+$ for $\lambda_j=
k_1\beta_1+\cdots+k_d\beta_d\in L^{\circ} /L, k_1, \cdots,
k_d\in \mathbb{C}$ and $2\lambda_j\in L$ is isomorphic to a direct sum
of
$$\otimes_{i=1}^{d}V_{{\mathbb
Z}\be_{i}+k_{i}\be_{i}}^{\epsilon_{i}}$$ with sign ${\epsilon_i=
\{\pm\}, 1\leq i\leq d}$ such that the number of $i$ with
$\epsilon_i= -$ is even.

(5) $V_{\lambda_j+L}^-$ for $\lambda_j=
k_1\beta_1+\cdots+k_d\beta_d, k_1, \cdots, k_d\in \mathbb{C}$ and
$2\lambda_j\in L$ is isomorphic to a direct sum of
$$\otimes_{i=1}^{d}V_{{\mathbb
Z}\be_{i}+k_{i}\be_{i}}^{\epsilon_{i}}$$ with sign ${\epsilon_i=
\{\pm\}, 1\leq i\leq d}$ such that the number of $i$ with
$\epsilon_i= -$ is odd.

(6)   $(V_L^{T_\chi})^+$ is isomorphic to a direct sum of
$$(V_{\mathbb{Z}\beta_1}^{T_{\chi_1}})^{\epsilon_1}\otimes\cdots \otimes
(V_{\mathbb{Z}\beta_d}^{T_{\chi_d}})^{\epsilon_d}$$ with signs
$\epsilon_i\in \{\pm\}, i = 1, \cdots, d$ such that the number of
$i$ with $\epsilon_i= -$ is even.

(7) $(V_L^{T_\chi})^-$ is isomorphic to a direct sum of
$$(V_{\mathbb{Z}\beta_1}^{T_{\chi_1}})^{\epsilon_1}\otimes\cdots \otimes
(V_{\mathbb{Z}\beta_d}^{T_{\chi_d}})^{\epsilon_d}$$ with signs
$\epsilon_i\in \{\pm\}, i=1, \cdots, d$ such that the number of $i$
with $\epsilon_i= -$ is odd.

\end{lem}
\pf (1)   By Lemma \ref{lc5.1}, we have $V_L\cong \otimes_{i=
1}^{d}{V_{\mathbb{Z}\beta_i}}$ and the corresponding $\theta$ is
changed to $\theta_1\otimes \cdots \otimes \theta_d$, where
$\theta_i$ is the restriction of $\theta$ to
$V_{\mathbb{Z}\beta_i}$, then  $V_L^+$ is isomorphic to
$(\otimes_{i= 1}^{d}{V_{\mathbb{Z}\beta_i}})^+$ as $\otimes_{i=
1}^{d}{V_{\mathbb{Z}\beta_i}^+}$-modules. The decomposition of
$(\otimes_{i= 1}^{d}{V_{\mathbb{Z}\beta_i}})^+$ into direct sum of irreducible
$\otimes_{i=1}^{d}{V_{\mathbb{Z}\beta_i}^+}$-modules is now obvious.

The proof of (2) is similar to that of (1). (3) is obvious. For (4),
note that
$$
V_{\lambda_j+L}=V_{\Z\beta_1+k_1\beta_1}\otimes \cdots \otimes
V_{\Z\beta_d+k_d\beta_d}.$$ Since $2\la_{j}\in L$ and
$\{\be_{1},\cdots,\be_{d}\}$ is a basis of $L$, it follows that
$2k_{i}\be_{i}\in {\mathbb Z}\be_{i}$ and $V_{{\mathbb
Z}\be_{i}+k_{i}\be_{i}}=V_{{\mathbb Z}\be_{i}-k_{i}\be_{i}}$, for
$i=1,2,\cdots,d$. The decomposition then follows easily. The proof
of (5) is similar.

Now we consider  the last two cases. Note that $(\al,\be)\in 2\Z$ for all
$\al,\be\in L.$ From the discussion given in Section 3 (before Theorem \ref{t32})
we see that
$$V_L^{T_\chi}=V_{\Z\be_1}^{T_{\chi_1}}\otimes\cdots \otimes V_{\Z\be_d}^{T_{\chi_d}}.$$
(6) and (7) follows immediately. \qed

By Lemma \ref{l4.1c} we have:
\begin{lem}\label{ll5.1} The lowest weights of irreducible $V_L^+$-modules are given by
$$
\begin{tabular}{|c|c|c|c|c|c|}
\hline

$V_L^+$ & $V_L^-$ & $V_{\lambda_i+L}$ & $V_{\mu_j+L}^{\pm}$ & $(V_L^{T_{\chi_i}})^+$  & $(V_L^{T_{\chi_j}})^-$ \\
\hline

0 & 1 & $\frac{\langle \lambda_i,\lambda_i\rangle }{2}$ & $\frac{\langle \mu_j,\mu_j\rangle }{2}$ & $\frac{d}{16}$ & $\frac{d+8}{16}$ \\
\hline
\end{tabular}
$$
where $2\la_{i}\notin L$, $2\mu_{j}\in L$ and $\mu_{j}\neq 0$.
\end{lem}
From Lemma \ref{ll5.1}, Theorem \ref{n3.7}, Proposition \ref{t3.3c}
and Theorem \ref{t5.10}, we get the following result:
\newtheorem{name3}{Theorem}[section]
\begin{lem}\label{l5.3c}
$Ext_{V_L^+}^1(N, M) = 0$ for the following
 irreducible $V_L^+$-module pairs $(M, N)$:
$$(M,N) = (M, M)\ ({\rm i.e.}\  M = N),  \  (V_{\lambda_j+L}^{\pm}, \  V_{\lambda_j+L}^{\mp}), \ \la_{j}\neq 0 ,$$
$$((V_L^{T_{\chi_i}})^{\mp}, (V_L^{T_{\chi_j}})^{\pm}), \ ((V_L^{T_{\chi_i}})^{\pm}, (V_L^{T_{\chi_j}})^{\pm}).$$
\end{lem}

Furthermore, we have:
\begin{lem}\label{l5.4c}
The extension groups $Ext_{V_L^+}^1(V_L^{\pm},V_L^{\mp})=0.$
\end{lem}
\pf  By Theorem \ref{n2.1} and Proposition \ref{n3.8}, we only need to prove $$Ext_{V_L^+}^1(V_L^-,V_L^+) = 0.$$
We consider an exact sequence $$0\rightarrow
V_L^-\rightarrow M\rightarrow V_L^+\rightarrow 0$$ for a weak
$V_L^+$-module $M$. By the rationality of
$\otimes_{i=1}^{d}{V_{\mathbb{Z}\beta_i}^+}$, there exists a
$\otimes_{i=1}^{d}{V_{\mathbb{Z}\beta_i}^+}$-submodule $M^1$ of $M$ such that
$M^1\cong V_L^+$ as $\otimes_{i=1}^{d}{V_{\mathbb{Z}\beta_i}^+}$-modules.
Since $\otimes_{i=1}^{d}{V_{\mathbb{Z}\beta_i}^+}$ and $V_L^+$ have the same Virasoro element, then there is a vector $u$ in $M^1$ such that
$L(-1)u= 0$ and $L(0)u=0$, this implies that $u$ generates a
$V_L^+$-submodule isomorphic to $V_L^+$, then we have $M\cong
V_L^-\bigoplus V_L^+$ and $Ext_{V_L^+}^1(V_L^-, V_L^+)=0$. \qed

\vskip 0.2cm
 We next prove  $Ext_{V_L^+}^1(M^2, M^1) = 0$ for the
remaining pairs $(M^1, M^2)$.
\begin{lem}\label{l5.5c}
The extension groups $Ext_{V_L^+}^1(M^2, M^1) = 0$ for the following
pairs $(M^1, M^2)$:
$$(M^1,M^2)= (V_{\lambda_i+L}, V_{\la_j+L}) \ \la_{i}\neq \la_{j}, \ (V_{\lambda_i+L},  V_{\la_j+L}^{\pm}),  \ \ (V_{\lambda_i+L}^{\pm}, V_{\la_j+L}),$$
$$ (V_{\lambda_i+L}^{\pm}, V_{\la_j+L}^{\pm}), \ \ (V_{\lambda_i+L}^{\pm}, V_{\la_j+L}^{\mp}), \ \la_{i}\neq \la_{j} ,$$
$$((V_L^{T_{\chi}})^{\pm}, V_{\la_{j}+L}^{\pm}), \ \  ((V_L^{T_{\chi}})^{\mp},  V_{\la_{j}+L}^{\pm}), \
(V_{\la_{j}+L}^{\mp}, (V_L^{T_{\chi}})^{\pm}), \ \
(V_{\la_{j}+L}^{\pm}, (V_L^{T_{\chi}})^{\pm}),$$
$$(V_{\lambda_i+L}, (V_L^{T_{\chi}})^{\pm}), \ \
((V_L^{T_{\chi}})^{\pm}, V_{\lambda_i+L}).$$
\end{lem}
\pf Let $(M^1,M^2)$ be one of the following pairs
$$((V_L^{T_{\chi}})^{\pm}, V_{\la_{j}+L}^{\pm}), \ \  ((V_L^{T_{\chi}})^{\mp},  V_{\la_{j}+L}^{\pm}), \
(V_{\la_{j}+L}^{\mp}, (V_L^{T_{\chi}})^{\pm}), \ \
(V_{\la_{j}+L}^{\pm}, (V_L^{T_{\chi}})^{\pm}),$$
$$(V_{\lambda_i+L}, (V_L^{T_{\chi}})^{\pm}), \ \
((V_L^{T_{\chi}})^{\pm}, V_{\lambda_i+L}).$$
 Let
$U=\otimes_{i=1}^{d}{V_{\mathbb{Z}\beta_i}^+}$. Then by Theorem
\ref{t2.6} and Theorem \ref{t3.3d}, $U$ is rational. It is obvious
that  $U$ has the same Virasoro element with $V$. Then by Lemma
\ref{la1}, it is enough to show that
$$I_{U}\left(\begin{tabular}{c}
$N^1$\\
$N$\  $N^2$
\end{tabular}\right)=0$$
for any irreducible $U$-submodules $N^1,N,N^2$ of $M^1,V_L^+,M^2,$
respectively. By (1)-(2) and (6)-(7) of Lemma \ref{l5.2c}, we know
that there is exactly one $N^i$ that has the form
$$(V_{\mathbb{Z}\beta_1}^{T_{\chi_1}})^{\epsilon_1}\otimes\cdots \otimes
(V_{\mathbb{Z}\beta_d}^{T_{\chi_d}})^{\epsilon_d}.$$ Then by Theorem
\ref{3.20} (2) and Theorem \ref{n2.2c}, the fusion rule of type
$$\left(\begin{tabular}{c}
$N^1$\\
$N$ $N^2$\\
\end{tabular}\right)$$
for the vertex operator algebra
$U=\otimes_{i=1}^{d}{V_{\mathbb{Z}\beta_i}^+}$ is zero.  Thus
$Ext_{V_L^+}^1(M^2, M^1) = 0$.

Consider the pair $(M^1, M^2) = (V_{\lambda_i+L},V_{\la_j+L})$,
where $\la_{i}=k_1\beta_1+\cdots +k_d\beta_d$,
$\la_{j}=l_1\beta_1+\cdots +l_d\beta_d$ such that
$\la_{i}\not=\la_{j}$ and $2\la_{i}, 2\la_{j}\notin L.$  Without
loss of generality, we may assume that $k_1\ne l_1.$ By (3) of Lemma
\ref{l5.2c}, we have
$$ V_{\lambda_i+L}=V_{\mathbb{Z}\beta_1+k_1\beta_1}
\otimes\cdots\otimes V_{\mathbb{Z}\beta_d+k_d\beta_d}, $$
$$V_{\lambda_j+L}= V_{\mathbb{Z}\beta_1+l_1\beta_1}
\otimes\cdots\otimes V_{\mathbb{Z}\beta_d+l_d\beta_d}. $$ Note that
$V_{\mathbb{Z}\beta_r+n_r\beta_r}$ is an irreducible
$V_{\mathbb{Z}\beta_r}$-module if
$2n_r\beta_r\not\in\mathbb{Z}\beta_r$ and
$V_{\mathbb{Z}\beta_r+n_r\beta_r}=V_{\mathbb{Z}\beta_r+n_r\beta_r}^{+}\bigoplus
V_{\mathbb{Z}\beta_r+n_r\beta_r}^{-}$ is a sum of two irreducible
$V_{\mathbb{Z}\beta_r}$-modules if $2n_{r}\be_{r}\in {\mathbb
Z}\be_{r}.$ Let $N^1,N,N^2$ be any irreducible $U$-submodules of
$M^1,V_L^+,M^2,$ respectively. It  follows immediately from Theorem
\ref{3.20} (2) and  Theorem \ref{n2.2c} that
$$I_{U}\left(\begin{tabular}{c}
$N^1$\\
$N$\  $N^2$
\end{tabular}\right)=0.$$ So by Lemma
\ref{la1}, $Ext_{V_L^+}^1{(M^2, M^1)}= 0$ in this case.

If $(M^1, M^2) = (V_{\lambda_i+L},V_{\la_j+L}^{\pm}), \
(V_{\lambda_i+L}^{\pm},V_{\la_j+L}),$ where $2\la_{i}\notin L$,
$2\la_{j}\in L$, by Proposition \ref{n3.8} and Theorem \ref{n2.1},
we only need to consider the pairs
$(M^1,M^2)=(V_{\lambda_j+L}^{\pm}, V_{\la_i+L})$. Let
$\lambda_j=k_1\beta_1+\cdots +k_d\beta_d$ and
$\la_i=l_1\beta_1+\cdots +l_d\beta_d$. Then $2k_{i}\beta_{i}\in
\Z\beta_{i}$, for $i=1,2,\cdots,d$. Since $2\la_{i}\notin L$,
it follows that there exists $1\leq s\leq d$ such that
$2l_{s}\beta_{s}\notin \Z\beta_{s}$.
By (3)-(5) of Lemma \ref{l5.2c},
$$V_{\lambda_i+L}= V_{\mathbb{Z}\beta_1+l_1\beta_1}
\otimes\cdots\otimes V_{\mathbb{Z}\beta_d+l_d\beta_d}, $$
 where
$V_{\mathbb{Z}\beta_r+l_r\beta_r}=V_{\mathbb{Z}\beta_r+l_r\beta_r}^{+}\bigoplus
V_{\mathbb{Z}\beta_r+l_i\beta_r}^{-}$ if
$2l_{r}\be_{r}\in {\mathbb Z}\be_{r}.$
Let $N^1,N,N^2$ be any
irreducible $U$-submodules of $M^1,V_L^+,M^2,$ respectively.
Then $N^{1}$ has the form:
$$\otimes_{i=1}^{d}V_{{\mathbb
Z}\be_{i}+k_{i}\be_{i}}^{\epsilon_{i}}.$$ We know from Theorem
\ref{3.20} (2) and  Theorem \ref{n2.2c} that
$$
I_{U}\left(\begin{tabular}{c}
$N^1$\\
$N$\  $N^2$
\end{tabular}\right)=0.$$ Again  by Lemma
\ref{la1}, we have $Ext_{V_L^+}^1{(M^2, M^1)}= 0$.

Finally we deal with the pairs $(M^1, M^2)$: $$ \
(V_{\lambda_i+L}^{\pm}, V_{\la_j+L}^{\pm}),\ (V_{\lambda_i+L}^{\pm},
V_{\la_j+L}^{\mp}),  \ \la_{i}\neq \la_{j}.$$  Since $\la_i\ne
\la_j,$ it follows that there exists $1\leq s\leq d$ such that
$k_s\ne l_s.$ Let $\lambda_i=k_1\beta_1+\cdots +k_d\beta_d$ and
$\la_j=l_1\beta_1+\cdots +l_d\beta_d$. Then $2k_{i}\beta_{i},
2l_{i}\beta_{i}\in \Z\beta_{i}$, for $i=1,2,\cdots, d$. Let
$N^1,N,N^2$ be any irreducible $U$-submodules of $M^1,V_L^+,M^2,$
respectively. By (3)-(4) of Lemma \ref{l5.2c}, $N^{1}$ and $N^{2}$
have the form:
$$V_{\mathbb{Z}\beta_1+k_1\beta_1}^{\epsilon_{i}} \otimes\cdots\otimes
V_{\mathbb{Z}\beta_d+k_d\beta_d}^{\epsilon_{i}}, $$ and
$$V_{\mathbb{Z}\beta_1+l_1\beta_1}^{\epsilon_{i}} \otimes\cdots\otimes
V_{\mathbb{Z}\beta_d+l_d\beta_d}^{\epsilon_{i}} $$ respectively. By
Theorem \ref{3.20} (2) and Theorem \ref{n2.2c},
$$
I_{U}\left(\begin{tabular}{c}
$N^1$\\
$N$\  $N^2$
\end{tabular}\right)=0.$$
This implies that $Ext_{V_L^+}^1{(M^2,
M^1)}= 0$. The
proof is complete. \qed

We are now in a position to state the main result of this section.

\begin{theorem}\label{n5.1}
Let $L$ be a positive definite even lattice with an orthogonal base,
then $V_L^+$ is rational.
\end{theorem}
\pf The theorem follows from Lammas \ref{l5.3c}-\ref{l5.5c} and
Theorem \ref{n2.2}.\qed

\section{Rationality of $V_L^+:$ general case}

We are now in a position to deal with $V_L^+$ for any positive definite even lattice
$L$ by using the rationality result obtained in Section 5.

Note that $L$ has a sublattice $L_1= \bigoplus_{i=
1}^d{\mathbb{Z}\beta_i}$ of the same rank, where $\{\beta_i, 1\leq i
\leq d\}$ is an orthogonal subset in $L$
 in the sense that
$\langle\beta_i, \beta_j\rangle= 0$, $i\neq j.$ This induces an
embedding
$$V_{L_1}^+= (V_{\bigoplus_{i=1}^d{\mathbb{Z}\beta_i}})^+\subset
V_L^+$$ of vertex operator algebras. Let $\{\gamma_s\}$ be a set of
representatives of $L/L_1$.  In Section 5, we prove that $V_{L_1}^+$
is rational. Then we can decompose all the irreducible
$V_L^{+}$-modules as a direct sum of irreducible
$V_{L_1}^+$-modules. Specifically, we have
\begin{lem}\label{l6.1}
(1) \  $V_L^+$ is isomorphic to a direct sum of
$$V_{\gamma_s+L_1}^+,(2\gamma_s\in
L_1), V_{\gamma_s+L_1}, (2\gamma_s\notin L_1)$$ (one of the  two
irreducible modules may not exist in the direct sum).

(2)  $V_L^-$ is isomorphic to a direct sum of
$$V_{\gamma_s+L_1}^-, (2\gamma_s\in
L_1), V_{\gamma_s+L_1}, (2\gamma_s\notin L_1)$$ (one of the  two
irreducible modules may not exist in the direct sum).

(3) \  $V_{\lambda_j+L}^+$ for $(2\lambda_j \in L)$ is isomorphic to
a direct sum of
$$V_{\lambda_j+\gamma_s+L_1}^+ \  (2(\la_{j}+\gamma_s)\in
L_1),\  V_{\lambda_j+\gamma_s+L_1} (2(\la_{j}+\gamma_s)\notin L_1)$$
(one of the  two irreducible modules may not exist in the direct
sum).

(4) \ $V_{\lambda_j+L}^-$ for $(2\lambda_j \in L)$ is isomorphic to
a direct sum of
$$V_{\lambda_j+\gamma_s+L_1}^- \  (2(\la_{j}+\gamma_s)\in
L_1), \ V_{\lambda_j+\gamma_s+L_1}\  (2(\la_{j}+\gamma_s)\notin
L_1)$$ (one of the  two irreducible modules may not exist in the
direct sum).

(5) \  $V_{\lambda_j+L}$ for $(2\lambda_j\notin L)$ is isomorphic to
a direct sum of
$$V_{\lambda_j+\gamma_s+L_1}.$$

(6)  $(V_L^{T_{\chi}})^{+}$  is isomorphic to a direct sum of
irreducible $V_{L_1}^{+}$-modules
$$(V_{L_1}^{T_{\chi_i}})^{+}$$ for some irreducible $\widehat{L_1}/K_1$-module ${T_{\chi_i}}$
with central character $\chi_i$ such that $\chi_{i}(\kappa)=-1$,
where $K_1=\{\theta(a)a^{-1}|a\in {\widehat{L_1}}\}$ and
$\widehat{L_1} = \{a\in \widehat{L}|\overline{a}\in L_1\}$.

(7) $(V_L^{T_{\chi}})^-$ is isomorphic to a direct sum of
$$(V_{L_1}^{T_{\chi_i}})^-$$ for some irreducible $\widehat{L_1}/K_1$-module ${T_{\chi_i}}$ with central
character $\chi_i$  such that $\chi_{i}(\kappa) = -1$.
\end{lem}
\pf (1)-(4) are obvious. (5) is immediate by noting that
$2(\lambda_j+\gamma_s)\notin L_1$ as $L_1$ is a sublattice of $L.$
Now we prove (6)-(7). Let $\widehat{L_1} = \{a\in
\widehat{L}|\overline{a}\in L_1\}$ and $K_1 = \{\theta(a)a^{-1}|a\in
{\widehat{L_1}}\}$.  Then $\widehat{L_1}$ is an abelian group
isomorphic to $L_1\times \langle\kappa\rangle$ and $K_1$ is
isomorphic to $2L_1.$ As a result, $\widehat{L_1}/K_1$ is isomorphic
to $L_1/2L_1\times \langle\kappa\rangle.$ Since $\widehat{L_1}/K_1$
is a subgroup of $\widehat{L}/K_1,$ $T_{\chi}$ is a direct sum of
one-dimensional irreducible $\widehat{L}/K_1$-modules $T_{\chi_i}$
such that $\kappa$ acts as $-1.$ Since ${\mathfrak
h}=\C\otimes_{\Z}L_{1}=\C\otimes_{\Z} L$, we see that
$$V_{L}^{T_{\chi}}=M(1)(\theta)\otimes T_{\chi}=\bigoplus_{i}M(1)(\theta)\otimes T_{\chi_i}=\bigoplus_{i}V_{L_1}^{T_{\chi_i}}.$$
That is, $V_{L}^{T_{\chi}}$ is a direct sum of irreducible $\theta$-twisted $V_{L_1}$-modules $V_{L_1}^{T_{\chi_i}}.$
(6)-(7) are evident now. \qed

By  Lemma \ref{ll5.1}, Theorem \ref{n3.7}, Proposition \ref{t3.3c}
and Theorem \ref{t5.10}, we have the following  result similar to
Lemma \ref{l5.3c}.
\newtheorem{name4}{Theorem}[section]
\begin{lem}\label{l6.1c}
The extension groups $Ext_{V_L^+}^1(N, M)=0$ for the following
irreducible $V_L^+$-module pairs $(M, N)$,
$$(M,M) \ ({\rm i.e} \ M = N), \ (V_{\lambda_j+L}^{\pm}, V_{\lambda_j+L}^{\mp}), \ \la_{j}\neq 0,$$
$$((V_L^{T_{\chi_i}})^{\mp}, (V_L^{T_{\chi_j}})^{\pm}), \ ((V_L^{T_{\chi_i}})^{\pm}, (V_L^{T_{\chi_j}})^{\pm}).$$
\end{lem}

Note that $V_{L_1}^+$ and $V_{L}^+$ have the same Virasoro element.
An analogue of Lemma \ref{l5.4c} with the same proof is the following:
\begin{lem}\label{l6.2c}
The extension groups $Ext_{V_L^+}^1(V_L^{\pm}, V_L^{\mp}) = 0$.
\end{lem}

For the remaining pairs $(M^{1},M^{2})$, we also have the following
result:
\begin{lem}\label{l6.3c}
 $Ext_{V_L^+}^1(M^{2}, M^{1})=0$ for the following
irreducible $V_L^+$-module pairs   $(M^1, M^2)$:
$$(M^1,M^2)= (V_{\lambda_i+L}, V_{\la_j+L}) \ \la_{i}\neq \la_{j}, \ (V_{\lambda_i+L},  V_{\la_j+L}^{\pm}),  \ \ (V_{\lambda_i+L}^{\pm}, V_{\la_j+L}),$$
$$ (V_{\lambda_i+L}^{\pm}, V_{\la_j+L}^{\pm}), \ \ (V_{\lambda_i+L}^{\pm}, V_{\la_j+L}^{\mp}), \ \la_{i}\neq \la_{j} ,$$
$$((V_L^{T_{\chi}})^{\pm}, V_{\la_{j}+L}^{\pm}), \ \  ((V_L^{T_{\chi}})^{\mp},  V_{\la_{j}+L}^{\pm}), \
(V_{\la_{j}+L}^{\mp}, (V_L^{T_{\chi}})^{\pm}), \ \
(V_{\la_{j}+L}^{\pm}, (V_L^{T_{\chi}})^{\pm}),$$
$$(V_{\lambda_i+L}, (V_L^{T_{\chi}})^{\pm}), \ \
((V_L^{T_{\chi}})^{\pm}, V_{\lambda_i+L}).$$
\end{lem}
\pf  Let $(M^1, M^2)$ be one of the following pairs:
$$((V_L^{T_{\chi}})^{\pm}, V_{\la_{j}+L}^{\pm}), \ \  ((V_L^{T_{\chi}})^{\mp},  V_{\la_{j}+L}^{\pm}), \
(V_{\la_{j}+L}^{\mp}, (V_L^{T_{\chi}})^{\pm}), \ \
(V_{\la_{j}+L}^{\pm}, (V_L^{T_{\chi}})^{\pm}),$$
$$(V_{\lambda_i+L}, (V_L^{T_{\chi}})^{\pm}), \ \
((V_L^{T_{\chi}})^{\pm}, V_{\lambda_i+L}).$$ By Theorem \ref{n5.1}
and Lemma \ref{la1}, it is enough to show that
$$I_{V_{L_{1}}^{+}}(\begin{tabular}{c}
$N^1$\\
$N$\  $N^2$
\end{tabular})=0$$
for any irreducible $V_{L_{1}}^{+}$-submodules $N^1,N,N^2$ of
$M^1,V_L^+,M^2,$ respectively.  By (6)-(7) of Lemma \ref{l6.1},
there is exactly one $N^i$ which has the form
$$(V_{L_{1}}^{T_{\chi_{i}}})^{+} \ \ or \ \
(V_{L_{1}}^{T_{\chi_{i}}})^{-}.$$ By Theorem \ref{t3.7}, we have
$$ I_{V_{L_{1}}^{+}}\left(\begin{tabular}{c}
$N^1$\\
$N$\  $N^2$
\end{tabular}\right)=0.$$

Now consider the pair $(M^1,M^2)=(V_{\lambda_i+L}, V_{\la_j+L})$
with $2\la_{i}\notin L, 2\la_{j}\notin L, \la_{i}\neq \la_{j}.$ By
the choice of $\{\la_{k}|k\in L^{\circ}/L\}$ (see Section 3), we
have $L+\lambda_i\neq L+\lambda_j$.  If  $L+\lambda_i=L-\lambda_j$,
then by  Lemma \ref{ll5.1} $V_{\lambda_i+L}$ and $V_{-\la_j+L}$ have
the same lowest weight. By Proposition \ref{t3.3c} and Theorem
\ref{t5.10}, $Ext_{V_L^+}^1(M^{2}, M^{1})=0.$ So we now assume that
$L+\la_i\neq L\pm \la_j.$ Let $N^1,N,N^2$ be any irreducible
$V_{L_{1}}^{+}$-submodules of $M^1,V_L^+,M^2,$ respectively. Then by
(5) of Lemma \ref{l6.1},
$$
N^{1}=V_{\lambda_i+\gamma_{s}+L_{1}}, \
N^{2}=V_{\la_j+\gamma_{r}+L_{1}},$$  for some
$\gamma_{r},\gamma_{s}$, where $2(\lambda_i+\gamma_s)\notin L_{1}$,
$2(\lambda_j+\gamma_r)\notin L_{1}$. By (1) of Lemma \ref{l6.1}, $N$
is of the form $V_{\gamma_{l}+L_{1}}^{+} (2\gamma_{l}\in L_{1})$ or
$V_{\gamma_{l}+L_{1}} (2\gamma_{l}\notin L_{1})$.
We claim that $(\gamma_{l}, \la_{i}+\gamma_{s}, \la_{j}+\gamma_{r})$
is not an admissible triple modulo $L_{1}$. Otherwise,
since $\gamma_{r}, \gamma_{l}, \gamma_{s}\in  L$ and $L_{1}\subseteq L$, this forces
$\la_{i}+\la_{j}\in L$   or $\la_{i}-\la_{j}\in L,$ a contradiction.
Since $(\gamma_{l}, \la_{i}+\gamma_{s}, \la_{j}+\gamma_{r})$
is not an admissible triple modulo $L_{1},$
Theorem \ref{t3.7} asserts that
$$I_{V_{L_{1}}^{+}}\left(\begin{tabular}{c}
$N^1$\\
$N$\  $N^2$
\end{tabular}\right)=0.$$

For the pairs $(M^1,M^2)=(V_{\lambda_i+L}$,  $V_{\la_j+L}^{\pm}),
(V_{\lambda_j+L}^{\pm}, V_{\la_i+L})$, where $2\la_{i}\notin L$,
$2\la_{j}\in L$, by Proposition \ref{n3.8} and Theorem \ref{n2.1},
we only need to consider the pairs
$(M^1,M^2)=(V_{\lambda_j+L}^{\pm}, V_{\la_i+L})$. Let $N^1,N,N^2$ be
any irreducible $U$-submodules of $M^1,V_L^+,M^2,$ respectively.
Then by (3)-(5) of Lemma \ref{l6.1}, $N^{2}=
V_{\la_i+\gamma_{r}+L_{1}}$, for some $\gamma_{r}$ and $N^{1}$ is
one of the following:
$$V_{\lambda_j+\gamma_s+L_1}^{\pm} \  (2(\la_{j}+\gamma_s)\in
L_1), \ V_{\lambda_j+\gamma_s+L_1} \  (2(\la_{j}+\gamma_s)\notin
L_1).$$ By (1) of Lemma \ref{l6.1}, $N$ is of the form
$V_{\gamma_{l}+L_{1}}^{+} (2\gamma_{l}\in L_{1})$ or
$V_{\gamma_{l}+L_{1}} (2\gamma_{l}\notin L_{1})$.
Note that $(\gamma_{l}, \la_{j}+\gamma_{s}, \la_{i}+\gamma_{r})$  is
not an admissible triple modulo $L_{1}$. Otherwise, either
$\la_i+\la_j$ or $\la_i-\la_j\in L.$ In either case we conclude
that $2\la_i\in L,$ a contradiction.
Again by Theorem
\ref{t3.7},
$$I_{V_{L_{1}}^{+}}\left(\begin{tabular}{c}
$N^1$\\
$N$\  $N^2$
\end{tabular}\right)=0.$$

The proof for the remaining pairs $(M^1,M^2)$:
 $$ (V_{\lambda_i+L}^{\pm}, V_{\la_j+L}^{\pm}), \
(V_{\lambda_i+L}^{\pm}, V_{\la_j+L}^{\mp})$$ is quite similar. We
omit it. \qed

\vskip 0.2cm
 By Lemmas \ref{l6.1c}-\ref{l6.3c} and Theorem \ref{n2.2},
together with Theorem \ref{n5.1}, we have the following main result
of the paper.

\begin{theorem}
Let $L$ be a positive definite even lattice. Then $V_{L}^{+}$ is
rational.
\end{theorem}

\end{document}